**Improved SQP and SLSQP Algorithms for Feasible Path-based Process Optimisation**

Yingjie Ma[1,2#*], Xi Gao[3,4#], Chao Liu[1] and Jie Li[1*]

[1]Centre for Process Integration, Department of Chemical Engineering, School of Engineering, The University of Manchester, Manchester, UK, M13 9PL

[2]Department of Chemical Engineering, Massachusetts Institute of Technology, Cambridge, MA, 02139

[3]School of Electronic and Information Engineering, Tongji University, Shanghai, China 201804

[4]School of Mechanical and Electrical Engineering, Jinggangshan University, Ji'an, Jiangxi, China 343009

**Abstract**

Feasible path algorithms have been widely used for process optimisation due to its good convergence. The sequential quadratic programming (SQP) algorithm is usually used to drive the feasible path algorithms towards optimality. However, existing SQP algorithms may suffer from inconsistent quadratic programming (QP) subproblems and numerical noise, especially for ill-conditioned process optimisation problems, leading to a suboptimal or infeasible solution. In this work, we propose an improved SQP algorithm (I-SQP) and an improved sequential least squares programming algorithm (I-SLSQP) that solves a least squares (LSQ) subproblem at each major iteration. A hybrid method through the combination of two existing relaxations is proposed to solve the inconsistent subproblems for better convergence and higher efficiency. We find that a certain part of the dual LSQ algorithm suffers from serious cancellation errors, resulting in an inaccurate search direction or no viable search direction generated. Therefore, the QP solver is used to solve LSQ subproblems in such a situation. The computational results indicates that I-SLSQP is more robust than fmincon in MATLAB, IPOPT, Py-SLSQP and I-SQP. It is also shown that I-SLSQP and Py-SLSQP is superior to I-SQP for ill-conditioned process optimisation problems, whilst I-SQP is more computationally efficient than I-SLSQP and Py-SLSQP for well-conditioned problems.

**Keywords**: Process optimisation; feasible path algorithm; numerical noise; ill conditioning; sequential least squares programming algorithm; sequential quadratic programming.

[#] Contributed equally to this work.

[*] To whom correspondence should be addressed. jie.li-2@manchester.ac.uk. Tel: +44 (0) 161 529 3084 (Jie Li). yingjma@mit.edu (Yingjie Ma).

# 1 Introduction

Optimisation is a powerful tool to design the best chemical process with the lowest cost or highest profit while satisfying the production requirements and restrictions at the same time (Biegler, 1993). To get optimisation results that can match the real-world production well, it is highly desirable to use rigorous unit operation models (Biegler, 2010). However, this often leads to strongly nonlinear, non-convex or even ill-conditioned large-scale nonlinear programming (NLP) problems, which are challenging to solve.

A typical NLP problem is provided as follows,

$$\min_{\mathbf{x} \in \mathbb{R}^n} f(\mathbf{x}) \qquad \text{(NLP)}$$

$$\text{s.t.} \quad \mathbf{h}(\mathbf{x}) = 0,$$

$$\mathbf{g}(\mathbf{x}) \geq 0,$$

where $\mathbf{x}$ is a vector of real variables with $n$ dimensions, and $f: \mathbb{R}^n \to \mathbb{R}$, $\mathbf{g}: \mathbb{R}^n \to \mathbb{R}^{m_I}$ and $\mathbf{h}: \mathbb{R}^n \to \mathbb{R}^{m_E}$ are functions that are at least twice continuously differentiable. $f$ is called the objective function, while $\mathbf{g}(\mathbf{x}) \geq 0$ and $\mathbf{h}(\mathbf{x}) = 0$ are inequality constraints and equality constraints, respectively.

There are mainly four categories of methods for solving process optimisation problems, including stochastic algorithms, surrogate model-based optimisation, simultaneous methods, and feasible/infeasible path methods. The former two algorithms (Ledezma-Martínez et al., 2018; Caballero & Grossmann, 2008) can exploit existing commercial process simulators more conveniently, but they are usually slow and cannot guarantee solution quality well. The simultaneous methods formulate a large and sparse NLP problem that consists of process model equations (equalities in the NLP problem), process requirements (equalities and inequalities) and design objective, and then solve the problem using an existing NLP solver. Such large and sparse NLP problems can be solved very efficiently by modern NLP algorithms, such as generalized reduced gradient (GRG) methods (Drud, 1994), interior point algorithms (Byrd et al., 1999), and sparse sequential quadratic programming (SQP) algorithms (Gill et al., 2002) developed in the last several decades if good initial points are available. However, finding a good initial point may take a significant amount of effort that could be prohibitive for practitioners and researchers (Kossack et al., 2006; Dowling & Biegler, 2015; Ma et al., 2021). If the algorithms do not converge, any physically meaningful solution cannot be obtained. Instead, the feasible/infeasible path algorithms are widely used to solve process optimisation problems due to their good convergence (Biegler, 2010), which decompose the entire problem into a small-scale NLP problem in the outer level and a large-scale process simulation problem in the inner level (Biegler et al., 1981; Biegler and Hughes, 1982). Hence, the feasible/infeasible path algorithms are the focus of the current work. The infeasible path algorithms differ from the feasible path algorithms in that the former put the equalities corresponding to all or part of recycle streams in the process flowsheet in the outer optimisation problem directly instead of in the inner simulation problems. The choice between feasible and infeasible path algorithms depends on the trade-off between solving a harder inner problem or facing a more challenging outer problem. Moreover, the feasible path algorithms generate a converged process simulation at each iteration, so the termination point can still be useful for chemical engineers even if the optimisation algorithms do not converge.

In both the feasible and infeasible path algorithms, the outer-level optimisation algorithm drives the inner-level simulations towards an optimal design iteratively, whilst the simulation provides the necessary information for the outer-level optimisation algorithm to determine a suitable step direction and length. Therefore, both simulation and optimisation algorithms are crucial for the feasible/infeasible path algorithms. For the inner-level simulation, the pseudo-transient continuation (PTC) modelling approach has been introduced to resolve the convergence issue of the equation-oriented simulation (Pattison and Baldea, 2014). We then combine the PTC simulation and the steady-state simulation to achieve a much higher efficiency (Ma et al., 2020a). Since the robust and efficient solution method for inner-level simulation problems has already existed, we will mainly focus on the feasible path algorithms. However, much more effort is required to develop a robust and efficient outer-level optimisation algorithm tailored for the feasible path algorithms, which will also benefit infeasible path algorithms.

In a feasible path algorithm, each function/gradient evaluation in the outer-level problem is conducted based on the results of a process simulation, which is usually quite expensive to evaluate. Therefore, most of its computational time is consumed in process simulations. Since the sequential quadratic programming (SQP) algorithm usually requires the least number of function evaluations (Powell, 1978b; Schittkowski, 1980), it is often used to solve the outer-level optimisation problem, although all the NLP algorithms that can be used in the simultaneous methods can also be used to drive the feasible path algorithms. SQP uses a quadratic programming (QP) problem to approximate problem (NLP) at an iterate $\mathbf{x}^k$ and generate the search direction $\mathbf{d}$ to find the next iterate $\mathbf{x}^{k+1}$, which is closer to the optimal solution $\mathbf{x}^*$ of problem (NLP). As a result, a sequence of iterates $\{\mathbf{x}^k\}$ is generated that presumably converges to $\mathbf{x}^*$ (Boggs and Tolle, 1995). This SQP method was first developed by (Wilson, 1963) and then modernized and popularized by a series of contributions (Han, 1976; Han, 1977) and (Powell, 1978a, b). In the Wilson-Han-Powell algorithm, the quasi-Newton method was used to approximate the Hessian matrix of the Lagrangian function in the QP subproblem to achieve superlinear local convergence. The line search and merit function were adopted to stabilize the algorithm and achieve global convergence.

When using SQP to drive the feasible path algorithm, one major issue is that the function and gradient values returned from process simulations contain numerical noise, as process simulations are regarded to be converged when some specified tolerance (e.g. $10^{-5}$) is satisfied (Moré and Wild, 2011). The numerical noise might be augmented significantly if the problem to be solved is ill-conditioned, which is not rare for real-world process optimisation problems (Biegler and Cuthrell, 1985). In such a situation, the search direction $\mathbf{d}$ generated from the QP subproblem could be inaccurate or even incorrect (e.g. ascent direction) (Schittkowski, 2011), resulting in many iterations required or even divergence. Although reducing the conditional number by appropriately scaling variables, constraints, and the objective function in problem (NLP) may help alleviate the issue, it is difficult to derive a general scaling method that can always improve the optimisation performance (Biegler and Cuthrell, 1985). Usually, it needs heuristics and trial and error to get good scaling factors (Ma et al., 2019). (Dai and Schittkowski, 2008) proposed a nonmonotone line search method that compared the current merit function value with the largest merit function value in the latest $t$ ($t \geq 2$) iterations, which could increase the possibility of finding a step length satisfying the line search condition. (Oztoprak et al., 2021) employed a relaxed line search strategy in the SQP algorithm to solve equality-constrained NLP problems with numerical noise. However, these methods have only been applied to small-scale examples (usually within 100 variables and constraints without any simulation-based implicit functions), such as those in Schittkowski (2008), instead of real-world process optimisation problems. Instead, in our previous contributions (Ma et al., 2020a; Ma et al., 2020b), we used a sequential least squares programming (SLSQP) algorithm (Kraft, 1988) to drive our hybrid steady-state and time-relaxation-based feasible path algorithms and successfully solved several challenging and ill-conditioned process optimisation problems. This indicates that SLSQP is quite promising for solving ill-conditioned problems whose reduced Hessian and/or Jacobian have large condition numbers (e.g. greater than $10^6$).

SLSQP, a variant of SQP, was first proposed by (Schittkowski, 1982). Instead of solving a QP subproblem at each major iteration, SLSQP solves an equivalent least squares (LSQ) subproblem to generate the search direction. (Schittkowski, 1982) reported that SLSQP needs a greater number of function evaluations than SQP due to less accurate descent directions generated according to the test results on a set of academic NLP problems. Until now, there are two SLSQP implementations (Schittkowski, 1982; Kraft, 1988), which both used a modified Powell's method to relax inconsistent QP subproblems (Powell, 1978b). However, the use of the modified Powell's method could not completely avoid premature termination (Tone, 1983). More seriously, the dual algorithm (Lawson and Hanson, 1995) used to solve LSQ subproblems in the existing SLSQP implementations may generate an ascending direction or an infeasible solution even though a feasible solution of the subproblem does exist.

With the above in mind, the existing SLSQP algorithm is improved and a robust SQP algorithm is developed concurrently in this work. Both the improved SLSQP and SQP algorithms use a hybrid relaxation method through the integration of the modified Powell's method (Powell, 1978b) and the Nowak's method (Nowak, 1988) to solve the inconsistent QP/LSQ subproblems. The former relaxation introduces one relaxation variable to relax all the constraints, leading to simpler QP/LSQ subproblems, while the Nowak's relaxation applies different relaxation variables for different constraints, providing larger flexibility. We analyse the dual LSQ solution algorithm and demonstrate the reason why it may incorrectly generate an infeasible or inaccurate solution. Therefore, in the improved SLSQP algorithm, when the dual LSQ solver reports an infeasible solution or generates an ascent/abnormal search direction even if the Hessian matrix has been reset as the identity matrix, the QP solver is activated to resolve the subproblem and generate a descent direction. Several challenging process optimisation problems are solved to validate the convergence and efficiency of the proposed algorithms in comparison to the fmincon solver in MATLAB (The Mathworks, 2023), the SLSQP solver (Kraft, 1988) in Python (Python Software Foundation, 2016) and the IPOPT solver (Wächter and Biegler, 2006).

## 2 Overview of SQP/SLSQP algorithms

### 2.1 Some notations and concepts

Throughout the work, $i$, $j$ and $k$ denote the index of variables, the index of constraints and the iteration number, respectively. $\mathcal{E}$ and $\mathcal{I}$ are the index set of equality and inequality constraints with dimensions of $m_E$ and $m_I$ respectively. All the vectors are column vectors. The symbol $\|\cdot\|$ denotes the 2-norm of a vector. $|\cdot|$ denotes the absolute value of a scalar or all the elements of a vector. $\nabla f$, $\nabla \mathbf{h}$ and $\nabla \mathbf{g}$ are the gradients of $f$, $\mathbf{h}$ and $\mathbf{g}$ respectively with $\nabla \mathbf{h} := [\nabla h_1, \nabla h_2, \ldots, \nabla h_{m_E}]$ and $\nabla \mathbf{g} := [\nabla g_1, \nabla g_2, \ldots, \nabla g_{m_I}]$. The active set of constraints is $\mathcal{A} := \mathcal{E} \cup \{j \in \mathcal{I} | g_j(x) = 0\}$. $L(\mathbf{x}, \boldsymbol{\lambda}, \boldsymbol{\mu}) = f(\mathbf{x}) + \boldsymbol{\lambda}^T \mathbf{h}(\mathbf{x}) - \boldsymbol{\mu}^T \mathbf{g}(\mathbf{x})$ is the Lagrangian function of problem (NLP), with $\boldsymbol{\lambda}$ and $\boldsymbol{\mu}$ being the Lagrange multipliers for equality and inequality constraints, respectively.

## 2.2 Basic SQP algorithm

The line search merit function based SQP algorithm for problem (NLP) solves the following quadratic programming problem (denoted as QP) to generate a search direction $\mathbf{d}$,

$$\min_{\mathbf{d} \in \mathbb{R}^n} \tfrac{1}{2}\mathbf{d}^T B^k \mathbf{d} + \nabla f(\mathbf{x}^k)^T \mathbf{d} \qquad \text{(QP)}$$

$$\nabla \mathbf{h}(\mathbf{x}^k)^T \mathbf{d} + \mathbf{h}(\mathbf{x}^k) = 0,$$

$$\nabla \mathbf{g}(\mathbf{x}^k)^T \mathbf{d} + \mathbf{g}(\mathbf{x}^k) \geq 0,$$

where $B^k$ denotes the approximate Hessian matrix of the Lagrangian function $L(\mathbf{x}^k, \boldsymbol{\lambda}^k, \boldsymbol{\mu}^k)$ with respective to $\mathbf{x}$. $\boldsymbol{\lambda}^k$ and $\boldsymbol{\mu}^k$ are the Lagrange multipliers of problem (QP). After solving problem (QP), the line search method is used to determine a suitable step length along the direction $\mathbf{d}$ generated. The iteration continues until specific convergence criteria are satisfied. A basic SQP algorithm is shown in **Algorithm 1**.

**Algorithm 1: A basic SQP algorithm**

Step 1:    $k \leftarrow 0$, given $\mathbf{x}^0$, $B^0$, $\boldsymbol{\rho}^0$, $\mathbf{v}^0$ and evaluate $f(\mathbf{x}^0)$, $\mathbf{g}(\mathbf{x}^0)$, $\mathbf{h}(\mathbf{x}^0)$, $\nabla f(\mathbf{x}^0)$, $\nabla \mathbf{g}(\mathbf{x}^0)$, $\nabla \mathbf{h}(\mathbf{x}^0)$;

Step 2:    solve problem (QP) to obtain the search direction $\mathbf{d}$ and Lagrange multipliers $\boldsymbol{\lambda}^k$, $\boldsymbol{\mu}^k$, and then proceed to the next step;

Step 3:    check the convergence criteria for problem (NLP). If they are satisfied, go to Step 7; otherwise, proceed to the next step;

Step 4:    update the penalty parameters $\boldsymbol{\rho}^k$, $\mathbf{v}^k$ using Eqs. (1-2) below and calculate the directional derivative $D\phi(\mathbf{x}^k, \mathbf{d}; \boldsymbol{\rho}^k, \mathbf{v}^k)$ using Eq. (4). Proceed to the next step;

Step 5:    conduct the line search with the merit function defined in Eq. (3) to get a step length $\alpha$ that satisfies the Armijio condition Eq. (5), set $\mathbf{x}^{k+1} \leftarrow \mathbf{x}^k + \alpha \mathbf{d}$, evaluate $f(\mathbf{x}^{k+1})$, $\mathbf{g}(\mathbf{x}^{k+1})$, $\mathbf{h}(\mathbf{x}^{k+1})$, and then proceed to the next step;

Step 6:     evaluate $\nabla f(\mathbf{x}^{k+1})$, $\nabla \mathbf{g}(\mathbf{x}^{k+1})$, $\nabla \mathbf{h}(\mathbf{x}^{k+1})$; update $B^{k+1}$ by damped BFGS formula in Eqs. (6-10); set $k \leftarrow k+1$, and then return to Step 2;

Step 7:     return $\mathbf{x}^k, \boldsymbol{\lambda}^k, \boldsymbol{\mu}^k, f(\mathbf{x}^k)$.

The penalty parameters are calculated using Eqs. (1-2),

$$\boldsymbol{\rho}^k = \max(|\boldsymbol{\lambda}^k|, \frac{\boldsymbol{\rho}^{k-1} + |\boldsymbol{\lambda}^k|}{2}), \tag{1}$$

$$\boldsymbol{\nu}^k = \max\left(|\boldsymbol{\mu}^k|, \frac{\boldsymbol{\nu}^{k-1} + |\boldsymbol{\mu}^k|}{2}\right), \tag{2}$$

where $\boldsymbol{\rho}^k$ and $\boldsymbol{\nu}^k$ are penalty parameters for the equality and inequality constraints, respectively. The L1 merit function $\phi(\mathbf{x}; \boldsymbol{\rho}^k, \boldsymbol{\nu}^k)$ is frequently used, which is defined by Eq. (3),

$$\phi(\mathbf{x}; \boldsymbol{\rho}^k, \boldsymbol{\nu}^k) = f(\mathbf{x}) + \sum_{j \in \mathcal{E}} \rho_j^k |h_j(\mathbf{x})| + \sum_{j \in \mathcal{I}} \nu_j^k g_j(\mathbf{x})^-, \tag{3}$$

where $g_j(\mathbf{x})^- := \max(0, -g_j(\mathbf{x}))$.

The directional derivative of the merit function along the direction $\mathbf{d}$ is:

$$D\phi(\mathbf{x}^k, \mathbf{d}; \boldsymbol{\rho}^k, \boldsymbol{\nu}^k) = \nabla f(\mathbf{x}^k)^T \mathbf{d} - \sum_{j \in \mathcal{E}} \rho_j^k |h_j(\mathbf{x}^k)| - \sum_{j \in \mathcal{I}} \nu_j^k g_j(\mathbf{x})^-, \tag{4}$$

The Armijio condition used in the line search is to guarantee the merit function has an enough decrease in each iteration to achieve the global convergence, as demonstrated in Eq. (5).

$$\phi(\mathbf{x}^k + \alpha \mathbf{d}; \boldsymbol{\rho}^k, \boldsymbol{\nu}^k) - \phi(\mathbf{x}^k; \boldsymbol{\rho}^k, \boldsymbol{\nu}^k) < \alpha \cdot \eta \cdot D\phi(\mathbf{x}^k, \mathbf{d}; \boldsymbol{\rho}^k, \boldsymbol{\nu}^k), \tag{5}$$

where $\eta \in (0, 0.5)$ is a constant. Here, we set it to be 0.1 throughout the work.

The following damped BFGS update (Nocedal and Wright, 2006; Powell, 1978c) is widely used as it demonstrates good performance and can guarantee positive definiteness of $B^k$ throughout optimisation.

$$B^{k+1} = B^k + \frac{\mathbf{r}^k (\mathbf{r}^k)^T}{(\mathbf{r}^k)^T \mathbf{s}^k} - \frac{B^k \mathbf{s}^k (\mathbf{s}^k)^T B^k}{(\mathbf{s}^k)^T B^k \mathbf{s}^k}, \tag{6}$$

where $\mathbf{s}^k = \mathbf{x}^{k+1} - \mathbf{x}^k$ and $\mathbf{r}^k = \theta \mathbf{y}^k + (1 - \theta) \mathbf{s}^k$. Here, $\mathbf{y}^k$ is the derivative change of the Lagrangian function between $\mathbf{x}^k$ and $\mathbf{x}^{k+1}$ defined in Eqs. (7-9), and the parameter $\theta$ is calculated by Eq. (10),

$$\mathbf{y}^k = \nabla_{\mathbf{x}} L(\mathbf{x}^{k+1}, \boldsymbol{\lambda}^k, \boldsymbol{\mu}^k) - \nabla_{\mathbf{x}} L(\mathbf{x}^k, \boldsymbol{\lambda}^k, \boldsymbol{\mu}^k), \tag{7}$$

$$\nabla_{\mathbf{x}} L(\mathbf{x}^{k+1}, \boldsymbol{\lambda}^k, \boldsymbol{\mu}^k) = \nabla f(\mathbf{x}^{k+1}) - \nabla \mathbf{h}(\mathbf{x}^{k+1}) \boldsymbol{\lambda}^k - \nabla \mathbf{g}(\mathbf{x}^{k+1}) \boldsymbol{\mu}^k, \tag{8}$$

$$\nabla_{\mathbf{x}} L(\mathbf{x}^k, \boldsymbol{\lambda}^k, \boldsymbol{\mu}^k) = \nabla f(\mathbf{x}^k) - \nabla \mathbf{h}(\mathbf{x}^k) \boldsymbol{\lambda}^k - \nabla \mathbf{g}(\mathbf{x}^k) \boldsymbol{\mu}^k, \tag{9}$$

$$\theta = \begin{cases} 1 & \text{if } (\mathbf{s}^k)^T \mathbf{y}^k \geq 0.2 (\mathbf{s}^k)^T B^k \mathbf{s}^k, \\ \dfrac{0.8 \cdot (\mathbf{s}^k)^T B^k \mathbf{s}^k}{(\mathbf{s}^k)^T B^k \mathbf{s}^k - (\mathbf{s}^k)^T y^k} & \text{if } (\mathbf{s}^k)^T \mathbf{y}^k < 0.2 (\mathbf{s}^k)^T B^k \mathbf{s}^k. \end{cases} \tag{10}$$

## 2.3 Reset the Hessian matrix

When solving some ill-conditioned optimisation problems with numerical noise, the following two scenarios may occur.

i) The line search may fail to generate a step length that satisfies the Armijo condition no matter how small the step length is;

ii) The search direction generated from problem (QP) may be an ascent direction for the merit function, which causes failure in optimisation if no action is taken.

To continue the optimisation, (Biegler and Hughes, 1985; Schittkowski, 2011) suggested to reset the Hessian matrix as the identity matrix and then resolve the QP subproblem to generate a new search direction. This strategy assumes that the problem is caused by the ill-conditioned Hessian matrix. However, if the numerical noise of function evaluations results in the unsatisfied Armijo condition, the line search may still fail even if the Hessian matrix is reset. In this case, the existing algorithm terminates the optimisation (Biegler and Hughes, 1985). Moreover, resetting a good Hessian approximation to identity leads to more optimisation iterations. Based on our extensive computational experience, however, we believe it is more robust and efficient to reset the Hessian matrix only when an ascent direction is generated, while accepting the last step length when reaching the maximum number of line searches even if the Armijo condition is not met. This strategy was used in the SLSQP code of (Kraft, 1988).

## 2.4 Convergence criteria

Theoretically, the Karush-Kuhn-Tucker (KKT) conditions (Nocedal and Wright, 2006) should be used as the convergence criteria. However, the KKT conditions cannot be satisfied within the required tolerance for many real-world problems due to its scale variance. Therefore, some other criteria are usually used in practice and often differ among works (Gill et al., 2019).

In this work, we use the same convergence criteria as those in the SLSQP code in Scipy (Kraft, 1988; Virtanen et al., 2020). There are two groups of criteria. The first group of criteria including Eqs. (11-13) are checked after obtaining the solution to the problem (QP).

$$acc_{inf} = \sum_{j\in\mathcal{E}} |h_j(\mathbf{x}^k)| + \sum_{j\in\mathcal{I}} g_j(\mathbf{x}^k)^- < tol, \tag{11}$$

$$acc_{opt} = |\nabla f(\mathbf{x}^k)^T \mathbf{d}| + |\boldsymbol{\lambda}^k|^T |\mathbf{h}(\mathbf{x}^k)| + |\boldsymbol{\mu}^k|^T \mathbf{g}(\mathbf{x}^k)^- < tol, \tag{12}$$

$$acc_{step} = \|\mathbf{d}\| < tol, \tag{13}$$

where $acc_{inf}$, $acc_{opt}$ and $acc_{step}$ represent the feasibility, optimality, and step length, respectively. $acc_{inf}$ is the summation of infeasibilities in all the constraints. $acc_{opt}$ indicates the decrease potential of the objective function and the weighted constraint infeasibility. $acc_{step}$ is the 2-norm of the descent direction. The solution is declared to be optimal if Eqs. (11-12) or Eqs. (11, 13) are satisfied.

The second group of convergence criteria involving Eqs. (14-16) are checked after the line search,

$$\widetilde{acc}_{inf} = \sum_{j\in\mathcal{E}} |h_j(\mathbf{x}^k + \alpha\mathbf{d})| + \sum_{j\in\mathcal{I}} g_j(\mathbf{x}^k + \alpha\mathbf{d})^- < \widetilde{tol}, \tag{14}$$

$$\widetilde{acc}_{opt} = |f(\mathbf{x}^k + \alpha\mathbf{d}) - f(\mathbf{x}^k)| < \widetilde{tol}, \tag{15}$$

$$acc_{step} = \|\mathbf{d}\| < \widetilde{tol}. \tag{16}$$

The optimisation is claimed to be successful if Eqs. (14-15) or Eqs. (14, 16) are satisfied. Normally, $\widetilde{tol} = tol$. However, when the search direction is ascent even if the Hessian matrix is the identity matrix or the Hessian matrix has been reset for certain times (e.g. $\bar{\iota}_{reset} = 5$ times), it is better to use a larger tolerance $\widetilde{tol} = \tau \cdot tol$, such as $\tau = 10$. This is because both scenarios usually indicate that the search region is close to the optimum and the influence of numerical noises in that region overwhelms the potential decrease in the merit function. Therefore, the strict tolerance $tol$ is quite difficult to achieve (Gill et al., 2019).

**Remark 1**: The returned solution is a feasible solution of the original NLP problem, if any group of convergence criteria is satisfied.

**Remark 2**: As discussed in (Gill et al., 2019), it is difficult to propose an optimality criterion that is generally suitable for all NLP problems to be solved, even those addressed using the same algorithm. Therefore, it is important to analyse the reasonability and optimality of the solution according to users' domain knowledge.

## 2.5 Basic SLSQP algorithm

The SLSQP algorithm basically follows the SQP algorithm except for steps 2 and 6. In step 2, instead of solving problem (QP), SLSQP solves the following linear constrained least squares subproblem (LSQ) to generate a descent direction **d**.

$$\min_{\mathbf{d}\in\mathbb{R}^n} \frac{1}{2}\|R^k\mathbf{d} - \mathbf{q}^k\|^2 \tag{LSQ}$$

$$\text{s.t. } [\nabla \mathbf{h}(\mathbf{x}^k)]^T \cdot \mathbf{d} + \mathbf{h}(\mathbf{x}^k) = 0,$$

$$[\nabla \mathbf{g}(\mathbf{x}^k)]^T \cdot \mathbf{d} + \mathbf{g}(\mathbf{x}^k) \geq 0,$$

where $R^k$ is an upper triangular matrix and $\mathbf{q}^k$ is a vector, which satisfy Eqs. (17) and (18) below respectively.

$$(R^k)^T R^k = B^k, \tag{17}$$

$$(R^k)^T \mathbf{q}^k = -\nabla f(\mathbf{x}^k). \tag{18}$$

After getting the $LDL^T$ factors ($L^k$ and $D^k$) of $B^k$, we can calculate the matrix $R^k$ by,

$$R^k = (D^k)^{\frac{1}{2}} \cdot (L^k)^T, \tag{19}$$

where $L^k$ is a lower triangular matrix with all the diagonal elements being 1, while $D^k$ is a diagonal matrix.

The other difference between SQP and SLSQP is related to the use of the BFGS formula in step 6. To improve computational accuracy and efficiency, the SLSQP algorithm updates $L^k$ and $D^k$ directly from the BFGS formula Eq. (20) instead of forming $B^k$ by Eq. (6) first and then factorizing it.

$$L^{k+1}D^{k+1}(L^{k+1})^T = L^k D^k (L^k)^T + \frac{\mathbf{r}^k(\mathbf{r}^k)^T}{(\mathbf{r}^k)^T \mathbf{s}^k} - \frac{L^k D^k (L^k)^T \mathbf{s}^k (\mathbf{s}^k)^T L^k D^k (L^k)^T}{(\mathbf{s}^k)^T L^k D^k (L^k)^T \mathbf{s}^k}, \tag{20}$$

The updating method is to apply the rank one modification algorithm from (Fletcher and Powell, 1974) twice since Eq. (20) is a rank two modification formula (Kraft, 1988). The details about the rank one updating method can be found in (Fletcher and Powell, 1974).

## 3 Improved SQP algorithm
### 3.1 Relaxations of QP subproblems

At some iterations of SQP, the QP subproblems may be infeasible even if the original problem (NLP) is feasible. Such QP subproblems are called inconsistent QP subproblems. To solve the inconsistent QP subproblems, the key idea is to solve a relaxation of the QP subproblem, which hence allows the SQP algorithm to continue.

There are two methods that are frequently used to develop such a relaxation of the QP subproblem. The first method is to introduce a slack variable $\xi \in [0,1]$ and construct the following relaxation (denoted as RQP1),

$$\min_{\mathbf{d}\in\mathbb{R}^n, \xi\in\mathbb{R}} \frac{1}{2}\mathbf{d}^T B^k \mathbf{d} + [\nabla f(\mathbf{x}^k)]^T \mathbf{d} + \frac{1}{2} M \cdot \xi^2 \quad \text{(RQP1)}$$

$$[\nabla \mathbf{h}(\mathbf{x}^k)]^T \mathbf{d} + \mathbf{h}(\mathbf{x}^k) - \xi \cdot \mathbf{h}(\mathbf{x}^k) = 0,$$

$$[\nabla \mathbf{g}(\mathbf{x}^k)]^T \mathbf{d} + \mathbf{g}(\mathbf{x}^k) - \xi \cdot C\mathbf{g}(\mathbf{x}^k) \geq 0,$$

$$0 \leq \xi \leq 1,$$

where $M$ is a constant to penalize the violation of the linear constraints. It is set to be $10^4$ throughout the work. $C$ is a diagonal matrix with diagonal elements defined as follows,

$$C_{j,j} = \begin{cases} 0 & \text{if } g_j(\mathbf{x}^k) > 0 \\ 1 & \text{if } g_j(\mathbf{x}^k) \leq 0 \end{cases}, \quad j \in \mathcal{J}. \quad (21)$$

Note that problem (RQP1) is always feasible because $\mathbf{d} = 0$ and $\xi = 1$ can satisfy its constraints trivially. However, when $\mathbf{d} = 0$, no progress in the optimisation will be achieved anymore. In such case, (Powell, 1978b) and (Biegler and Hughes, 1985) claimed the original problem (NLP) to be infeasible directly. However, such a conclusion may be incorrect as shown in (Tone, 1983).

To address this issue, another relaxation of the QP subproblems from Nowak (1988) is used, which is denoted as RQP2.

$$\min_{\substack{\mathbf{d}\in\mathbb{R}^n, \mathbf{s}\in\mathbb{R}^{m_E}, \\ \mathbf{t}\in\mathbb{R}^{m_E}, \mathbf{v}\in\mathbb{R}^{m_I}}} \frac{1}{2}\mathbf{d}^T B^k \mathbf{d} + \nabla f(\mathbf{x}^k)^T \mathbf{d} + \frac{1}{2} M'(\mathbf{s}^T\mathbf{s} + \mathbf{t}^T\mathbf{t} + \mathbf{v}^T\mathbf{v}) + \mathbf{w}_1^T(\mathbf{s} + \mathbf{t}) + \mathbf{w}_2^T \mathbf{v} \quad \text{(RQP2)}$$

$$[\nabla \mathbf{h}(\mathbf{x}^k)]^T \mathbf{d} + \mathbf{h}(\mathbf{x}^k) = \mathbf{s} - \mathbf{t},$$

$$[\nabla \mathbf{g}(\mathbf{x}^k)]^T \mathbf{d} + \mathbf{g}(\mathbf{x}^k) \geq -\mathbf{v},$$

$$\mathbf{s} \geq 0, \mathbf{t} \geq 0, \mathbf{v} \geq 0$$

where $\mathbf{s}$, $\mathbf{t}$ and $\mathbf{v}$ are three new vectors of variables introduced for relaxation of the constraints. $\mathbf{w}_1$ and $\mathbf{w}_2$ are two constant vectors used for penalizing constraint violations, and $M'$ is a constant for the same purpose.

The problem (RQP2) is similar to that of (Tone, 1983), but a second-order term was added in the objective function to ensure it is strictly convex (Nowak, 1988). Problem (RQP2) has a larger feasible region than problem (RQP1), allowing the SQP algorithm to have higher tendency to find a feasible solution (Tone, 1983). However, problem (RQP1) involves $2m_E + m_I - 1$ less variables than problem (RQP2), so the former can be solved more efficiently. Furthermore, more constraints are activated in problem (RQP2), which more often leads to the ill-conditioned constraint Jacobian. The (RQP2) problem with an ill-conditioned Jacobian is often difficult to solve. Both issues motivate the development of the following hybrid relaxation strategy in the following section.

## 3.2 Integration of the two relaxations of the QP subproblems

We propose an integration strategy to combine the advantages of the above two relaxations to solve an inconsistent QP subproblem during SQP. Once the QP subproblem is deemed infeasible, we first solve problem (RQP1) as problem (RQP1) has much fewer variables and is usually easier to solve. Once $\mathbf{d} = 0$ and $\xi = 1$ are generated from problem (RQP1), it indicates that no new iterate that leads to an improvement in the merit function could be found. As discussed before, it may be incorrect to assert that the original problem (NLP) is infeasible. We then resort to solving problem (RQP2). Otherwise, we examine the value of $\xi$ obtained from problem (RQP1). $\xi$ may be close to 1 (e.g. $\xi > \bar{\xi} = 0.99$) even if $\|\mathbf{d}\| > 0$, indicating that problem (RQP1) struggles in revising the inconsistent QP subproblem as large infeasibility exists. If this happens for more than a few times (e.g. $\bar{n} = 10$ times) consecutively, it means it is almost impossible to escape from the infeasible region of problem (QP) through solving problem (RQP1). We must resort to solving problem (RQP2).

When solving problem (RQP2), if the Jacobian matrix of the active constraints has a very large condition number $\kappa_A$ (e.g. $\kappa_A \geq \bar{\kappa} = 10^{30}$), which is usually considered to be singular, an incorrect solution or no solution may be generated, leading to slow optimisation progress or premature termination. If this occurs for more than a certain number of consecutive iterations (e.g. $\bar{n} = 10$ iterations), the relaxation is switched from (RQP2) back to (RQP1).

Finally, it should be highlighted that we solve the problem (QP) first in all cases. Problem (RQP1) or (RQP2) is solved only when problem (QP) is infeasible. The complete solution strategy using the hybrid relaxation for solving QP subproblems is shown in Fig. 1 and described below:

**Algorithm 2: Enhanced QP solution strategy using the hybrid relaxation**

Step 1:  Given $\bar{\xi}, \bar{n}, \bar{\kappa}$, the QP relaxation indicator $n_{rex}$ (1 for relaxation strategy 1, and 2 for relaxation strategy 2), number of continuous iterations ($n_\xi$) with $\xi \geq \bar{\xi}$, number of continuous iterations ($n_{ill}$) with $\kappa_A > \bar{\kappa}$;

Step 2:  Solve problem (QP). If it converges, go to step 13; otherwise, let $n_{rex} = 1$ and proceed to the next step;

Step 3  If $n_{rex} = 1$, solve problem (RQP1), set $n_{ill} \leftarrow 0$, and proceed to the next step; otherwise, solve problem (RQP2), set $n_\xi \leftarrow 0$, and go to Step 9;

Step 4  If problem (RQP1) converges, go to the next step; otherwise, go to step 13;

Step 5  If $\mathbf{d} = 0$, go to Step 8; otherwise, go to the next step;

Step 6  If $\xi < \bar{\xi}$, set $n_\xi \leftarrow 0$, go to Step 13; otherwise, go to the next step;

Step 7  If $n_\xi \geq \bar{n}$, go to Step 8; otherwise, set $n_\xi \leftarrow n_\xi + 1$, go to Step 13;

Step 8  Set $n_{rex} \leftarrow 2$, and then go back to Step 3;

Step 9  If problem (RQP2) converges, proceed to the next step; otherwise, go to step 13;

Step 10  If $\kappa_A < \bar{\kappa}$, set $n_{ill} \leftarrow 0$, go to Step 13; otherwise, go to the next step;

Step 11  If $n_{ill} \geq \bar{n}$, go to Step 12; otherwise, set $n_{ill} \leftarrow n_{ill} + 1$, go to Step 13;

Step 12  Set $n_{rex} \leftarrow 1$, and then go back to Step 3.

Step 13  Return (to the main SQP framework).

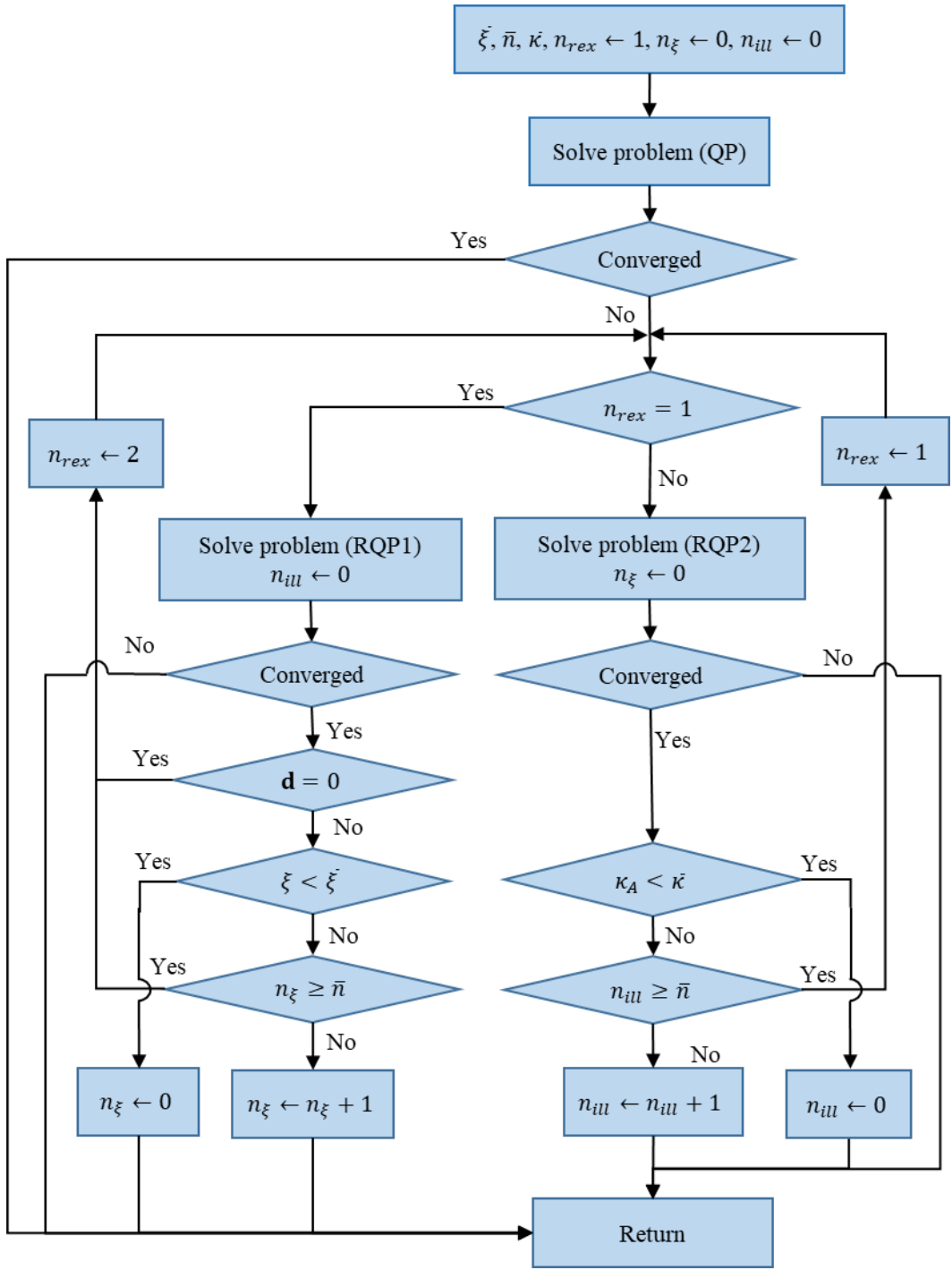

Figure 1 Flowchart of the enhanced QP solution procedure

The improved SQP algorithm with the hybrid relaxation strategy is illustrated in Fig. 2 and described in detail as follows.

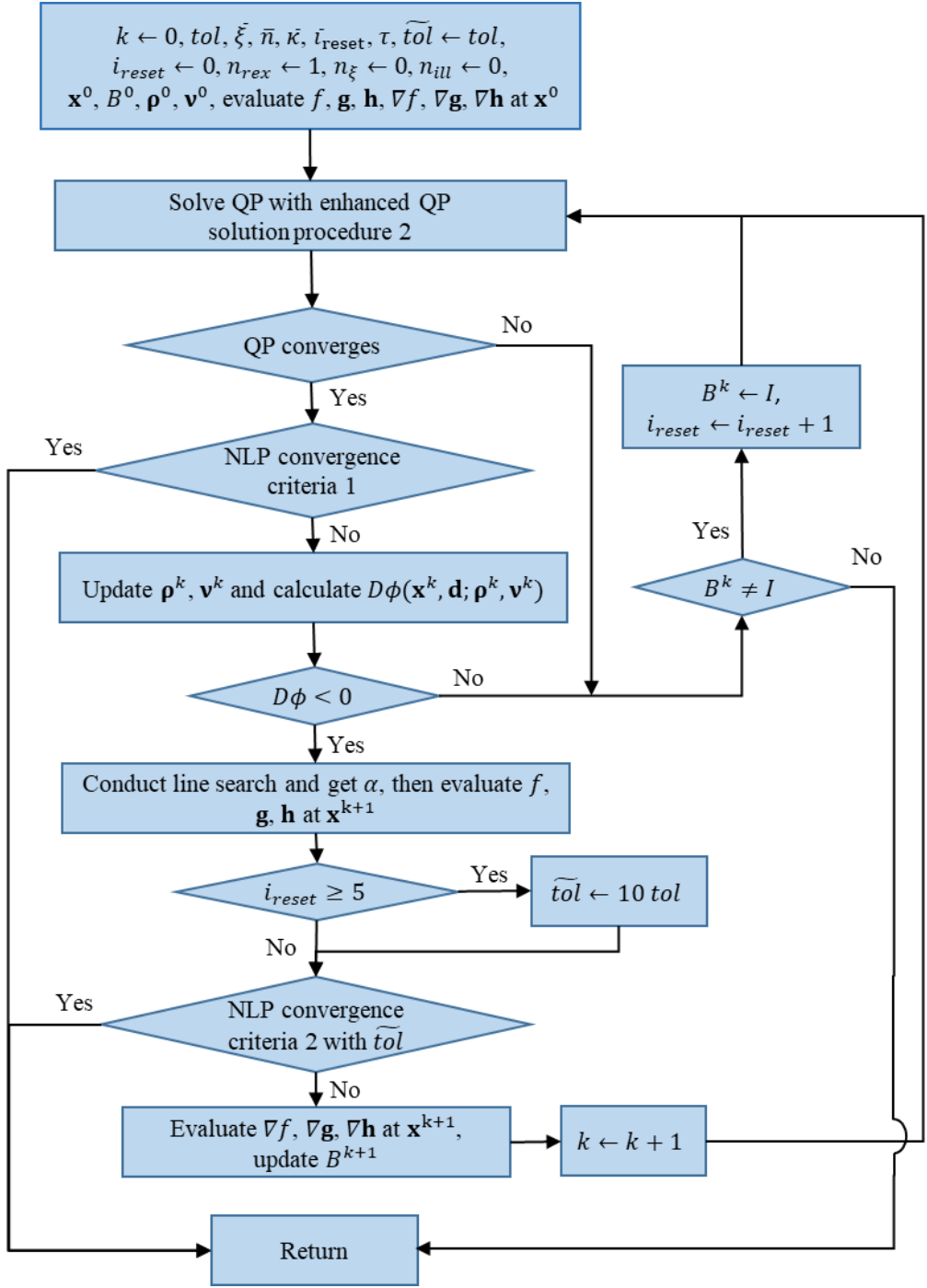

Figure 2 Flowchart of the improved SQP algorithm

**Algorithm 3: Improved SQP (I-SQP)**

Step 1: $k \leftarrow 0$, $tol$, $\bar{\xi}$, $\bar{n}$, $\bar{\kappa}$, $\bar{\iota}_{reset}$, $\tau$, $\widetilde{tol} \leftarrow tol$, $i_{reset} \leftarrow 0$, $n_{rex} \leftarrow 1$, $n_\xi \leftarrow 0$, $n_{ill} \leftarrow 0$, $\mathbf{x}^0$, $B^0$, $\boldsymbol{\rho}^0$, $\mathbf{v}^0$, evaluate $f(\mathbf{x}^0)$, $\mathbf{g}(\mathbf{x}^0)$, $\mathbf{h}(\mathbf{x}^0)$ $\nabla f(\mathbf{x}^0)$, $\nabla \mathbf{g}(\mathbf{x}^0)$, $\nabla \mathbf{h}(\mathbf{x}^0)$;

Step 2: solve the QP subproblem with the enhanced QP solution strategy (i.e. **Algorithm 2**) to get the search direction **d** and Lagrange multipliers $\boldsymbol{\lambda}^k$, $\boldsymbol{\mu}^k$. If no feasible solution is found, go to Step 6; otherwise, proceed to the next step;

Step 3: check the first group of convergence criteria for the original NLP problem. If it converges, go to Step 10; otherwise, proceed to the next step;

Step 4 update the penalty parameters $\boldsymbol{\rho}^k$, $\boldsymbol{v}^k$ using Eqs. (1-2) and calculate the directional derivative $D\phi(\mathbf{x}^k, \mathbf{d}; \boldsymbol{\rho}^k, \boldsymbol{v}^k)$ by Eq. (4). If $D\phi(\mathbf{x}^k, \mathbf{d}; \boldsymbol{\rho}^k, \boldsymbol{v}^k) \geq 0$, go to Step 6; otherwise, proceed to the next step;

Step 5 conduct the line search with the merit function defined in Eq. (3) to get a step length $\alpha$ satisfying Armijio condition Eq. (5), set $\mathbf{x}^{k+1} \leftarrow \mathbf{x}^k + \alpha \mathbf{d}$, evaluate $f(\mathbf{x}^{k+1})$, $\mathbf{g}(\mathbf{x}^{k+1})$, $\mathbf{h}(\mathbf{x}^{k+1})$, and then go to Step 7;

Step 6 if $B^k \neq \mathbf{I}$, $B^k \leftarrow \mathbf{I}$, $i_{reset} \leftarrow i_{reset} + 1$, and then return to Step 2; otherwise, go to Step 10;

Step 7 if $i_{reset} \geq \bar{\iota}_{reset}$, $\widetilde{tol} \leftarrow \tau \cdot tol$. Proceed to the next step;

Step 8 check the second group of convergence criteria for the original NLP problem. If it converges, go to Step 10; otherwise, proceed to the next step;

Step 9 evaluate $\nabla f(\mathbf{x}^{k+1})$, $\nabla \mathbf{g}(\mathbf{x}^{k+1})$, $\nabla \mathbf{h}(\mathbf{x}^{k+1})$, update $\mathbf{B}^{k+1}$ by the damped BFGS formula in Eqs. (6-10); set $k \leftarrow k + 1$, and then go back to Step 2;

Step 10 return $\mathbf{x}^k, \boldsymbol{\lambda}^k, \boldsymbol{\mu}^k, f(\mathbf{x}^k)$.

## 4  Improved SLSQP algorithm

As discussed before, the SLSQP algorithm solves a linear constrained least squares subproblem (LSQ) to generate a descent direction **d**, instead of solving a QP subproblem.

### 4.1  Inconsistent LSQ subproblems

Similar to SQP, inconsistent subproblems may be encountered during the iterations of SLSQP. The following two relaxations of the LSQ subproblem (denoted as RLSQ1 and RLSQ2 respectively) are used to tackle the inconsistent LSQ subproblems. The (RLSQ1) is a LSQ version of modified Powell's relaxation problem (Powell, 1978b), and it is implemented in the existing SLSQP codes (Schittkowski, 1982; Kraft, 1988). However, the (RLSQ1) may terminate prematurely at a point with $\mathbf{d} = 0$, so (RLSQ2) is required in such case, which is the counterpart of (RQP2). The use of (RLSQ2) is one advantage of the proposed SLSQP algorithm over existing ones.

$$\min_{\mathbf{d}\in\mathbb{R}^n, \xi\in\mathbb{R}} \frac{1}{2}\left\|\begin{bmatrix} R^k & \\ & M^{\frac{1}{2}} \end{bmatrix}\begin{bmatrix} \mathbf{d} \\ \xi \end{bmatrix} - \begin{bmatrix} \mathbf{q}^k \\ 0 \end{bmatrix}\right\|^2 \quad \text{(RLSQ1)}$$

$$s.t. \quad \nabla \mathbf{h}(\mathbf{x}^k)^T \cdot \mathbf{d} + \mathbf{h}(\mathbf{x}^k) - \xi \cdot \mathbf{h}(\mathbf{x}^k) = 0,$$

$$\nabla \mathbf{g}(\mathbf{x}^k)^T \cdot \mathbf{d} + \mathbf{g}(\mathbf{x}^k) - \xi \cdot \mathbf{c}^T \cdot \mathbf{g}(\mathbf{x}^k) \geq 0,$$

$$0 \leq \xi \leq 1.$$

$$\min_{\substack{\mathbf{d}\in\mathbb{R}^n, \mathbf{s}\in\mathbb{R}^{m_E}, \\ \mathbf{t}\in\mathbb{R}^{m_E}, \mathbf{v}\in\mathbb{R}^{m_I}}} \frac{1}{2}\left\|\begin{bmatrix} R^k & & & \\ & \left[M'^{\frac{1}{2}}\right]_{m_E} & & \\ & & \left[M'^{\frac{1}{2}}\right]_{m_E} & \\ & & & \left[M'^{\frac{1}{2}}\right]_{m_I} \end{bmatrix}\begin{bmatrix} \mathbf{d} \\ \mathbf{s} \\ \mathbf{t} \\ \mathbf{v} \end{bmatrix} - \begin{bmatrix} \mathbf{q}^k \\ -M'^{-\frac{1}{2}}\mathbf{w}_1 \\ -M'^{-\frac{1}{2}}\mathbf{w}_1 \\ -M'^{-\frac{1}{2}}\mathbf{w}_2 \end{bmatrix}\right\|^2 \quad \text{(RLSQ2)}$$

$$s.t. \quad \nabla \mathbf{h}(\mathbf{x}^k)^T \cdot \mathbf{d} + \mathbf{h}(\mathbf{x}^k) = \mathbf{s} - \mathbf{t},$$

$$\nabla \mathbf{g}(\mathbf{x}^k)^T \cdot \mathbf{d} + \mathbf{g}(\mathbf{x}^k) \geq -\mathbf{v},$$

$$\mathbf{s} \geq 0, \mathbf{t} \geq 0, \mathbf{v} \geq 0.$$

where $[M']_n$ is an $n$-by-$n$ diagonal matrix with all diagonal elements being $M'$. Like the enhanced QP solution strategy, we also propose an enhanced LSQ solution strategy through integration of the two relaxations RLSQ1 and RLSQ2 (i.e. **Algorithm 4**) as follows,

**Algorithm 4: Enhanced LSQ solution strategy using hybrid relaxations**

Step 1: Given $\bar{\xi}, \bar{n}, \bar{\kappa}$, the LSQ relaxation indicator $n_{rex}$ (1 for relaxation strategy 1, and 2 for relaxation strategy 2), number of continuous iterations ($n_\xi$) with $\xi \geq \bar{n}$, number of continuous iterations ($n_{\text{ill}}$) with $\kappa_A > \bar{\kappa}$;

Steps 2-12: The same as Steps 2-12 in **Algorithm 2** but substituting problems (QP), (RQP1), and (RQP2) with problems (LSQ), (RLSQ1), (RLSQ2) respectively;

Step 13: Return.

## 4.2 Numerical errors of the dual LSQ algorithm

There are several methods for solving LSQ subproblems. One obvious approach is to convert the LSQ subproblem to a QP subproblem, which, however, makes the SLSQP algorithm meaningless. The second method is the active-set null-space method, which solves an unconstrained LSQ problem in the null space of the linear constraints (Stoer, 1971). This method is similar to the frequently used active-set QP solution algorithm (Gill et al., 1984) and is expected to make SLSQP have similar performance to SQP. Another method is a dual algorithm proposed by (Lawson and Hanson, 1995), which is used in the SLSQP codes of (Schittkowski, 1982) and (Kraft, 1988). In this dual algorithm, the Householder Transformations (Golub and Van Loan, 2013) are used to eliminate equality constraints in problem (LSQ) and derive a linear least squares problem with only inequalities (LSI). In a further step, by substituting the objective function with a new variable $z$, the LSI problem is converted into a least distance problem (LDP).

$$\min_{\mathbf{z}\in\mathbb{R}^{n-m_E}} \frac{1}{2}\|\mathbf{z}\|^2 \qquad \text{(LDP)}$$

$$s.t. \quad \tilde{G}\mathbf{z} + \tilde{\mathbf{g}} \geq 0,$$

where $\tilde{G}$ is a $m_I$-by-$(n - m_E)$ matrix, and $\tilde{\mathbf{g}}$ is a $m_I$ dimensional vector. Finally, the dual counterpart of problem (LDP) is constructed as follows, which is called the nonnegative least squares (NNLS) problem (Lawson and Hanson, 1995),

$$\min_{\mathbf{u}\in\mathbb{R}^{m_I}} \frac{1}{2}\|A\mathbf{u} - b\|^2 \qquad \text{(NNLS)}$$

$$s.t. \; \mathbf{u} \geq 0.$$

Here, $A = \begin{bmatrix} \tilde{G}^T \\ -\tilde{\mathbf{g}}^T \end{bmatrix}$ and $\mathbf{b} = [\overbrace{0, \dots, 0}^{n'}, 1]^T$, where $n' = n - m_E$.

The dual algorithm proposed to solve problem (LSQ) suffers from numerical issues. First, the transformation from problem (LSI) to problem (LDP) is potentially unstable. As a result, infeasibilities in solving problem (NNLS) may be caused for some ill-conditioned problems even if the original problem (LSQ) is feasible, as shown in (Haskell and Hanson, 1981).

Second, a series of backward calculations are required to recover the solution $\mathbf{d}$ of the original problem (LSQ) from the solution $\mathbf{u}$ of the problem (NNLS). Numerical errors introduced when getting $\mathbf{z}$ from $\mathbf{u}$ by using the following Eqs. (22-23) may generate a $\mathbf{z}$ violating some of the linear constraints in problem (LDP), which consequently causes the search direction $\mathbf{d}$ to violate the linear constraints of the original problem (LSQ).

$$\mathbf{r} = A\mathbf{u} - \mathbf{b}, \qquad (22)$$

$$z_i = -\frac{r_i}{r_{n'+1}}, \qquad \forall i = 1,2,\dots n' \qquad (23)$$

where **r** is a vector of the residuals from solving problem (NNLS). Note that Eq. (22) potentially suffers from serious cancellation errors and the errors in **r** will be propagated into **z** as shown in **Appendix A**. The relationship between the computed solution $\mathbf{z} = [z_i]$ and the true solution $\mathbf{z}^* = [z_i^*]$ of the problem (LDP) are shown in Eq. (24).

$$z_i = z_i^* \frac{(1 \pm \epsilon_i^a)}{\left(1 \pm \frac{2\epsilon_{n'+1}^a}{r_{n'+1}^*}\right)}, i = 1,2,\dots,n', \text{when } r_{n'+1}^* \ll 1 \qquad (24)$$

where $\boldsymbol{\epsilon}^a$ are the errors in $A\mathbf{u}$, and $\mathbf{r}^*$ is the true solution of the problem (NNLS). According to Eq. (24), when the residual is quite small, e.g., $r_{n'+1}^* \leq 2\epsilon_{n'+1}^a$, the relative error of **z** may be nonsense; when the residual is a bit bigger, e.g. $r_{n'+1}^* \approx 20\epsilon_{n'+1}^a$, the relative error of **z** may be around 10%; when the residual is large enough, e.g. $r_{n'+1}^* \geq 200\epsilon_{n'+1}^a$, the relative error would be smaller than 1%. Here, $\boldsymbol{\epsilon}^a$ mainly includes the errors of the transformation process from problem (LSQ) to (NNLS) and the solution errors of problem (NNLS), so they're much larger than the machine precision (around $10^{-16}$ for double precision arithmetic operation). On the other hand, we observed $10^{-14}$ or even $10^{-17}$ for $r_{n'+1}$ during iterations. This means inaccurate or even wrong solutions might be generated from the LSQ algorithm during SLSQP iterations.

However, it seems that some extent of inaccuracy and instability of the dual algorithm can counteract the effect of the ill conditioning of the NLP problems, making the SLSQP algorithm perform evidently different from the SQP algorithm. But overly inaccurate LSQ solutions can cause slow progress or even premature termination of the optimisation. According to our observation, the following three issues may arise when applying the dual algorithm:

(1) The search direction with an abnormally large norm, which can be much larger (e.g. $\tau_d = 10$ times larger) than the norm of all the previous search directions. It usually causes the simulation to diverge during the line search, even if the PTC simulation is applied. As a result, a tiny step length is accepted. This slows down the optimisation.

(2) The ascent direction. The ascent direction is often caused by an ill-conditioned Hessian matrix, thus we can reset $L^k$ and $D^k$ as the identity matrix and then resolve the problem (LSQ). However, sometimes an ascent direction is generated even when $L^k$ and $D^k$ are identity matrices, which is due to the ill-conditioned Jacobian matrix of the active constraints. This can cause premature termination of the SLSQP algorithm.

(3) Failure to solve the LSQ subproblem even if it is feasible. Here, failure means both the original LSQ and the relaxed LSQ are reported to be infeasible. This could also cause the premature termination of the SLSQP algorithm.

To resolve the above issues, we first reset $L^k$ and $D^k$ as the identity matrix and then resolve the LSQ problem. If this does not work, a QP solver [e.g. the active set algorithm in Gurobi (Gurobi Optimisation, 2022)] is used to solve problem (LSQ).

The improved SLSQP algorithm is shown in Fig. 3 and described in detail as follows:

**Algorithm 5: Improved SLSQP algorithm (I-SLSQP)**

Step 1: $k \leftarrow 0$, $solver \leftarrow lsq$, $tol$, $\bar{\xi}$, $\bar{n}$, $\bar{\kappa}$, $\bar{\iota}_{reset}$, $\tau$, $\bar{k}$, $\tau_d$, $\widetilde{tol} \leftarrow tol$, $i_{reset} \leftarrow 0$, the maximum norm of descent direction $d_{max} \leftarrow 0$, $n_{rex} \leftarrow 1$, $n_\xi \leftarrow 0$, $n_{ill} \leftarrow 0$, $\mathbf{x}^0$, $L^0$, $D^0$, $\boldsymbol{\rho}^0$, $\mathbf{v}^0$, evaluate $R^0$, $\mathbf{q}^0$, $f(\mathbf{x}^0)$, $\mathbf{g}(\mathbf{x}^0)$, $\mathbf{h}(\mathbf{x}^0)$ $\nabla f(\mathbf{x}^0)$, $\nabla \mathbf{g}(\mathbf{x}^0)$, $\nabla \mathbf{h}(\mathbf{x}^0)$;

Step 2: solve the LSQ subproblem to get search direction $\mathbf{d}$ and Lagrange multipliers $\boldsymbol{\lambda}_k$, $\boldsymbol{\mu}_k$. When $solver = lsq$, solve the problem with the enhanced LSQ solution procedure 4; when $solver = qp$, solve the problem with the enhanced QP solution procedure 2. If no feasible solution is found, go to Step 7; otherwise, proceed to the next step;

Step 3: check the first group of convergence criteria for the NLP problem. If it converges, go to Step 12; otherwise, proceed to the next step;

Step 4: if $solver = lsq$ and $k \geq \bar{k}$ and $\|\mathbf{d}\| > \tau_d \cdot d_{max}$, go to step 7;

Step 5: update the penalty parameters $\boldsymbol{\rho}^k$, $\mathbf{v}^k$ using Eqs. (1-2) and calculate the directional derivative $D\phi(\mathbf{x}^k, \mathbf{d}; \boldsymbol{\rho}^k, \mathbf{v}^k)$ by Eq. (4). If $D\phi(\mathbf{x}^k, \mathbf{d}; \boldsymbol{\rho}^k, \mathbf{v}^k) \geq 0$, go to Step 7; otherwise, proceed to the next step;

Step 6: conduct the line search with the merit function defined in Eq. (3) to identify a step length $\alpha$ satisfying Armijio condition Eq. (5), set $\mathbf{x}^{k+1} \leftarrow \mathbf{x}^k + \alpha \mathbf{d}$, evaluate $f(\mathbf{x}^{k+1})$, $\mathbf{g}(\mathbf{x}^{k+1})$, $\mathbf{h}(\mathbf{x}^{k+1})$, and then go to Step 8;

Step 7: if $L^k \neq I$, $L^k \leftarrow I$, $D^k \leftarrow I$, $i_{reset} \leftarrow i_{reset} + 1$, and then return to Step 2; otherwise, go to Step 11;

Step 8: if $i_{reset} \geq \bar{\iota}_{reset}$, $\widetilde{tol} \leftarrow \tau \cdot tol$;

Step 9: check the second group of convergence criteria defined for the NLP problem; If it converges, go to Step 12;

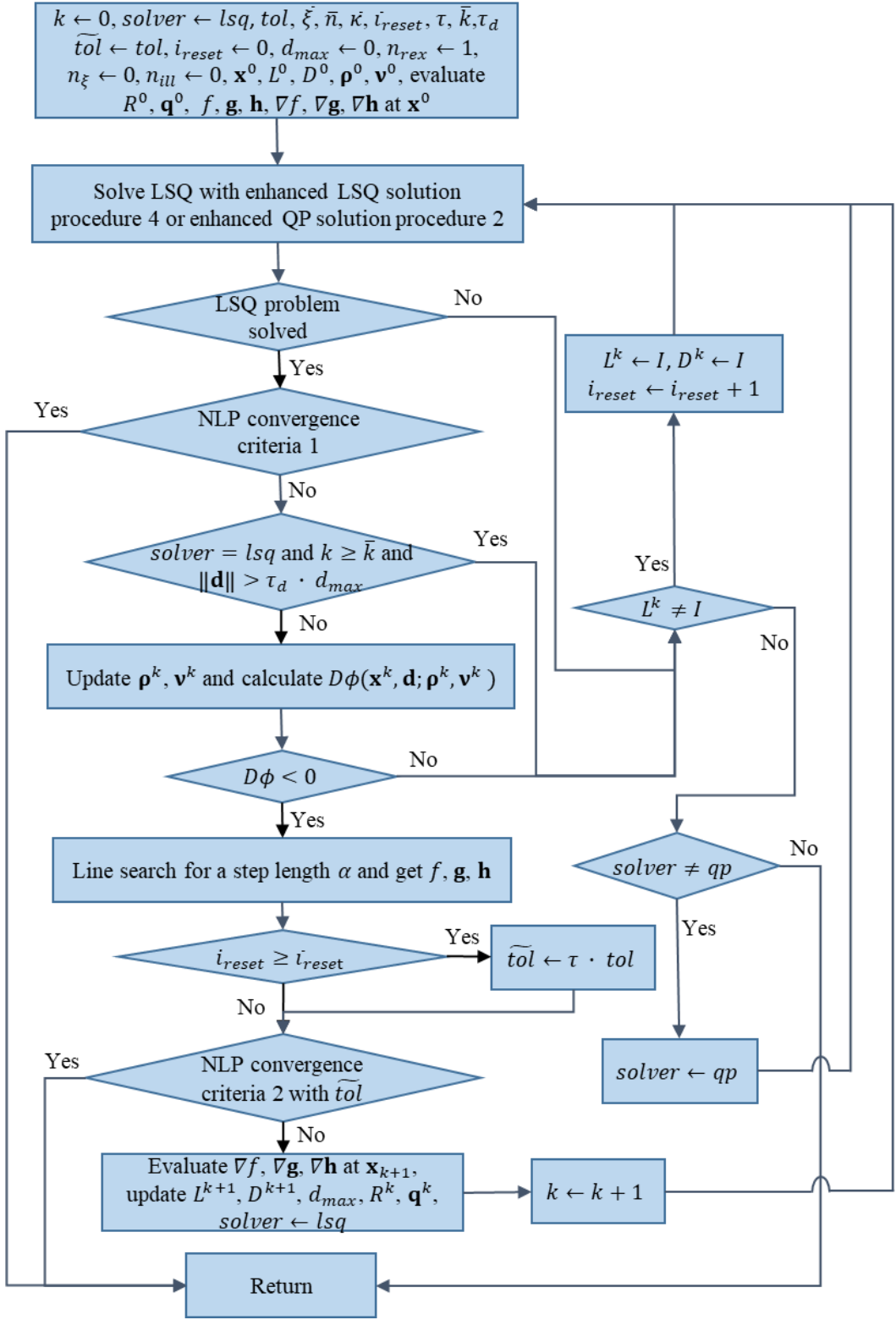

Figure 3 Flowchart of the improved SLSQP algorithm

Step 10: evaluate $\nabla f(\mathbf{x}^{k+1})$, $\nabla \mathbf{g}(\mathbf{x}^{k+1})$ $\nabla \mathbf{h}(\mathbf{x}^{k+1})$, update $L^{k+1}$ and $D^{k+1}$ by Eq. (20), and then get $R^k$ and $\mathbf{q}^k$ by Eqs (12-13); $d_{max} \leftarrow \max\{\|\mathbf{d}^{k'}\|, k' = 1,2,\ldots k\}$, $solver \leftarrow lsq$; set $k \leftarrow k+1$, and then go back to Step 2;

Step 11: If $solver \neq qp$, $solver \leftarrow qp$, and then go back to Step 2; otherwise, proceed to the next step;

Step 12 return $\mathbf{x}^k, \boldsymbol{\lambda}^k, \boldsymbol{\mu}^k, f(\mathbf{x}^k)$.

Note that in Step 4 the norm of the search direction is checked after $\bar{k}$ iterations (e.g. $\bar{k} = 5$ iterations) so that a reasonable value for $d_{max}$ is obtained after these iterations.

## 5 Computational studies

The hybrid steady-state and time-relaxation-based feasible path algorithm proposed in our previous work (Ma et al., 2020a) is used for optimisation. Process simulation is conducted in an equation-oriented environment such as Aspen Custom Modeler (Aspen Technology Inc., 2015), which provides the first-order derivatives based on algorithmic differentiation. The optimisation tolerance is 1×10$^{-5}$. The parameters used in the developed SQP and SLSQP algorithms are $\bar{\xi} = 0.99$, $\bar{n} = 10$, $\bar{\kappa} = 10^{30}$, $\bar{\iota}_{reset} = 5$ and $\tau = 10$, $\bar{k} = 5$, $\tau_d = 10$. Five algorithms are used to drive the feasible path algorithm respectively, including:

- Py-SLSQP: the SLSQP solver in Scipy 1.5.2, which is a wrapper of the SLSQP Fortran code written by (Kraft, 1988) with some improvements by the Scipy community (Virtanen et al., 2020). In our implementation, we restart the optimisation from the termination point at most 10 times when it terminates prematurely. This allows us to solve more problems.

- I-SQP: our improved SQP algorithm.

- I-SLSQP: our improved SLSQP algorithm.

- fmincon: SQP algorithm in the fmincon solver from Matlab 2023b (The Mathworks, 2023) is used here.

- IPOPT 3.14.13: interior point algorithm (Wächter and Biegler, 2006). The BFGS method is used to approximate the Hessian matrix. The linear solver MUMPS 5.2.1 is used. We did not use more efficient linear solvers for IPOPT because we realize that the computational time required to solve subproblems is usually negligible compared to that required for simulations in the feasible path algorithms.

Solution quality and computational efficiency are two important indicators that can be used to evaluate the performance of different algorithms. The former is evaluated by comparing optimal objective function values obtained from different algorithms. For evaluation of computational efficiency, as the QP subproblems can be solved very fast and nearly all the computational time is spent in process simulation, we use the total process simulation time (denoted as $t_{sim}$) as an index. It should be noted that we do not use the total number of function evaluations (denoted as $n_f$) as the primary evaluation criterion of computational efficiency. This is because a smaller number of function evaluations may still lead to a large computational cost if the time-consuming PTC simulations are used in the hybrid steady-state and time-relaxation based optimisation method.

We use the Morales profiles (Morales and Nocedal, 2011) to visualise the results from different algorithms for easy comparison. The Morales profiles of the following quantities obtained from two alternative algorithms (e.g. algorithm 1 vs. algorithm 2) are demonstrated.

$$\gamma_f = \log_2(\frac{f_1}{f_2}) \text{ and } \gamma_t = \log_2(\frac{t_1}{t_2}) \tag{25}$$

where $f_1$ and $f_2$ are the optimal objective function value from algorithm 1 (e.g. Py-SLSQP) and algorithm 2 (e.g. I-SLSQP) respectively. $t_1$ and $t_2$ denote the total computational time required by algorithm 1 (e.g. Py-SLSQP) and algorithm 2 (e.g. I-SLSQP), respectively. The negative $\gamma_f$ or $\gamma_t$ indicates that algorithm 1 is superior to algorithm 2, while a positive value shows that algorithm 2 is better. The values of $\gamma_f$ and $\gamma_t$ are then ranked individually in ascending order. The areas of the two half spaces can be used to evaluate the performance of the two algorithms. For visualisation of comparing the performance of different algorithms, the name of an algorithm is presented in the Morales profiles to denote the algorithm in the corresponding half space.

We evaluate the proposed I-SQP and I-SLSQP algorithms by solving seven large-scale process optimisation problems from the literature (Ma et al., 2020a; Ma et al., 2020b), which are challenging to solve. These seven problems cover a variety of intensified chemical processes including:

(1) a heat integrated pressure-swing distillation (PSD) process;

(2) a dividing-wall column with known and fixed pressure (DWC);

(3) a dividing-wall column with unknown pressure to be optimized (DWCP);

(4) a dividing-wall column-intensified reaction-separation-recycle process for production

of dimethyl ether (DME);

(5) an extractive dividing-wall column for separation of the acetone-chloroform mixture using dimethylsulfoxide (DMSO) as the solvent (EDWC-AC);

(6) an extractive dividing-wall column for separation of the ethanol-water mixture using ethylene glycol (EG) as solvent (EDWC-EG);

(7) a heat pump-assisted extractive dividing-wall column for separation of the ethanol-water mixture using EG as solvent (HPEDWC).

All these chemical processes are modelled using Aspen Custom Modeler V8.8. To solve each problem, six different initial points are generated to initialise the optimisation algorithms. These initial points differ by their initial bypass efficiencies of all stages in the distillation column, which are 0.1, 0.3, 0.5, 0.7, 0.9, and 1.0, respectively. The initial values of other decision variables are given in Appendix C. As a result, a total of 42 (7 × 6) problem instances are generated for evaluation. All instances solved by Py-SLSQP, I-SQP and I-SLSQP are executed on a desktop with a 3.20 GHz Intel® Core™ i7-8700 processor and 16 GB of RAM running Windows 10 64-bit operating system. All instances solved by fmincon are processed on a desktop with a 2.9 GHz Intel® Core™ i7-10700 processor and 32 GB of RAM running Windows 10 64-bit operating system. All instances solved by IPOPT are processed on a laptop with a 2.3 GHz Intel® Core™ i7-12700H processor and 16 GB of RAM running Windows 11 64-bit operating system.

The model statistics for all seven problems is provided in Table 1. The total number of function evaluations $n_f$, total simulation time $t_{sim}$ and optimum $f^*$ for each problem instance are presented in Tables 2 and 3.

Table 1 Model statistics for all seven problems

| Item | PSD | DWC | DWCP | DME | EDWC-AC | EDWC-EW | HPEDWC |
|---|---|---|---|---|---|---|---|
| $n_v$ | 6,411 | 21,352 | 21,352 | 13,661 | 10,570 | 15,805 | 16,117 |
| $n_{ind}$ | 62 | 186 | 187 | 122 | 90 | 125 | 128 |
| $n_{eq}$ | 5,727 | 18,970 | 18,970 | 12,150 | 9,399 | 14,038 | 14,295 |
| $n_{ieq}$ | 28 | 16 | 16 | 26 | 16 | 20 | 38 |

$n_v$: number of variables; $n_{ind}$: number of independent variables; $n_{eq}$: number of equalities; $n_{ieq}$: number of inequalities.

Table 2 Performance comparison of different algorithms/solvers for problems PSD, DWC and DWCP from six initial points

| Problem | Algorithm | No. of function evaluations | Time of function evaluations (s) | Optimal objective values |
|---|---|---|---|---|
| PSD (M\$ year$^{-1}$) | Py-SLSQP | 75/79/130/89/104/131 | 173/53/111/96/91/161 | 1.021/1.022/1.023/1.022/1.021/1.021 |
| | I-SQP | 54/160/77/76/90/82 | 193/193/127/74/98/114 | 1.021/1.022/1.022/1.022/1.021/1.023 |
| | I-SLSQP | 68/75/119/101/131/116 | 168/58/105/101/129/159 | 1.021/1.022/1.023/1.022/1.021/1.021 |
| | fmincon | 286/153/55/5/37/34 | 210/123/1294/na/na/146 | 1.005/1.005/inf/inf/inf/1.005 |
| | IPOPT | 505/538/228/256/422/475 | 386/248/132/156/307/314 | 1.047/1.046/1.047/1.047/1.005/1.005 |
| DWC (M\$ year$^{-1}$) | Py-SLSQP | 184/186/125/93/107/187 | 188/195/191/85/104/172 | 1.520/1.514/1.526/1.518/1.515/1.514 |
| | I-SQP | 69/69/182/64/66/72 | 84/88/348/64/65/72 | 1.517/1.515/1.515/1.515/1.515/1.515 |
| | I-SLSQP | 139/321/302/99/140/189 | 174/330/310/94/136/164 | 1.514/1.516/1.519/1.515/1.516/1.516 |
| | fmincon | 145/162/575/186/159/148 | 198/165/558/128/81/120 | 1.514/1.514/1.517/1.520/1.517/1.518 |
| | IPOPT | 1773/1667/1780/1761/1685/1548 | 1936/1724/1686/1848/1667/1590 | Inf/inf/inf/inf/inf/inf |
| DWCP (M\$ year$^{-1}$) | Py-SLSQP | 419/388/520/582/422/332 | 645/455/577/846/388/400 | 1.321/1.322/1.322/1.334/1.322/1.338 |
| | I-SQP | 717/163/239/204/370/457 | 933/186/231/220/322/408 | 1.587/1.326/1.331/1.323/1.325/1.334 |
| | I-SLSQP | 2478/269/607/1139/608/500 | 3890/381/922/1460/877/913 | 1.338/1.325/1.321/1.338/1.328/1.336 |
| | fmincon | 422/165/65/102/159/67 | 349/305/74/101/171/84 | Inf/inf/inf/inf/inf/inf |
| | IPOPT | 4540/4142/4458/4100/4075/3945 | 3543/3161/3649/2417/3470/2993 | inf/inf/inf/inf/inf/inf |

inf: infeasible solution; na: not available.

Table 3 Performance comparison of different solvers/algorithms for problems DME, EDWC-AC, EDWC-EW, and HPEDWC from six initial points

| Problem | Algorithm | No. of function evaluations | Time of function evaluations (s) | Optimal objective values |
| --- | --- | --- | --- | --- |
| DME $\left(\frac{TAC-123944}{10}\right)$ | Py-SLSQP | 608/819/1718/787/1306/988 | 1163/1231/7636/1867/3893/1481 | 1.740/1.766/1.766/1.768/1.761/1.766 |
| | I-SQP | 1819/2761/454/3234/631/810 | 1870/4343/1649/3845/670/1472 | 3.560/1.844/1.860/2.583/1.869/1.853 |
| | I-SLSQP | 310/951/741/1351/575/1071 | 538/1804/2454/3531/1827/4048 | 1.764/1.753/1.766/1.767/1.821/1.767 |
| | fmincon | 566/2001/11/570/1832/906 | na/1727/na/na/2800/551 | 1.915/inf/inf/inf/77.68/inf |
| | IPOPT | 2381/1919/2538/2987/617/1426 | 3719/2482/6243/6342/1469/na | Inf/inf/inf/inf/inf/inf |
| EDWC-AC ($10^5$ \$ year$^{-1}$) | Py-SLSQP | 249/185/164/234/278/394 | 855/953/516/918/1130/1326 | 6.081/6.105/6.191/6.109/7.962/6.077 |
| | I-SQP | 506/222/342/457/322/650 | 2485/234/734/969/898/2998 | 7.124/6.136/6.140/6.101/6.759/6.228 |
| | I-SLSQP | 335/205/303/168/256/356 | 1249/711/708/686/1413/1867 | 6.078/6.105/6.107/6.109/6.099/6.132 |
| | fmincon | 549/270/183/679/14/263 | 378/343/657/1357/na/1367 | 6.173/inf/inf/16.505/inf/inf |
| | IPOPT | 2613/2509/2328/2985/2594/1720 | 1653/3966/1134/4061/1648/953 | inf/inf/inf/inf/inf/inf |
| EDWC-EW (M\$ year$^{-1}$) | Py-SLSQP | 271/994/inf/487/211/778 | 801/3564/inf/1762/1528/3276 | 5.459/5.389/inf/5.390/5.393/5.381 |
| | I-SQP | 364/696/inf/327/784/600 | 1487/3472/inf/1185/3785/1974 | 5.386/5.407/inf/12.252/6.951/5.460 |
| | I-SLSQP | 353/657/374/223/361/649 | 974/1319/2165/967/3381/3538 | 5.402/5.382/5.435/5.541/5.415/5.383 |
| | fmincon | 2001/212/23/1/1/19 | 4405/1838/na/na/na/na | inf/inf/inf/inf/inf/inf |
| | IPOPT | 2161/2492/2352/1743/227/2963 | 3990/5627/6199/6694/na/6557 | Inf/inf/inf/inf/inf/inf |
| HPEDWC (M\$ year$^{-1}$) | Py-SLSQP | 881/978/1660/inf/415/814 | 2074/2821/5123/inf/1705/2638 | 4.826/4.927/4.830/inf/5.025/4.849 |
| | I-SQP | 642/562/606/667/318/973 | 2727/1635/1267/1981/1208/2780 | 5.542/5.867/4.939/5.019/5.134/4.985 |
| | I-SLSQP | 810/409/1152/806/545/442 | 1702/1316/2853/2435/1355/2432 | 4.833/4.836/4.828/4.829/4.839/5.039 |
| | fmincon | 1/393/764/1/1/15 | na/na/na/na/na/79 | Inf/inf/inf/inf/inf/inf |
| | IPOPT | 578/2550/2882/2948/2802/2522 | 4331/3450/3191/5304/4764/4526 | Inf/inf/inf/inf/inf/inf |

inf: infeasible solution; na: not available.

## 5.1  fmincon and IPOPT vs. PySLSQP, I-SQP and I-SLSQP

As seen from Tables 2 and 3, both fmincon and IPOPT have serious convergence issue when solving the given process optimisation problems. Specifically, fmincon can optimize the PSD problem and DWC problem from at least three different initial points, but it struggles in solving the other five problems. IPOPT is even worse, which could only solve the PSD problem. Instead, Py-SLSQP and I-SQP can solve most of the problems, while I-SLSQP can solve all problems. fmincon and IPOPT also need more function evaluations than the other three algorithms for the converged instances. Especially, IPOPT is the slowest among the algorithms in terms of the number of function evaluations, which indicates the interior point algorithm may not be suitable to be used in the feasible path algorithms, although it may be advantageous for the simultaneous optimization methods due to its computationally tractable subproblems (Biegler, 2010).

## 5.2  Py-SLSQP vs. I-SLSQP

We first compare the performance of Py-SLSQP and I-SLSQP. The comparative results are provided in Tables 2 and 3. As shown in Table 3, Py-SLSQP fails to solve the EDWC-EW problem from the third initial point and the HPEDWC problem from the fourth initial point due to the positive directional derivatives obtained at some intermediate iterations. However, I-SLSQP can solve all problems from all given initial points (i.e. all 42 problem instances).

The Morales profiles of $\gamma_f$ and $\gamma_t$ for the 40 problem instances with converged solutions from both Py-SLSQP and I-SLSQP are illustrated in Figure 4. Note that the two problem instances that Py-SLSQP fails to solve are not included. As shown in Fig. 4a, Py-SLSQP and I-SLSQP generate similar optimal solutions for most instances due to $\gamma_f$ being nearly zero from instances 7 to 34. The highest $\gamma_f$ is obtained in instance 40, which corresponds to solving the problem EDWC-AC from the fifth initial point. The corresponding TAC from I-SLSQP and Py-SLSQP is $6.099 \times 10^5$ \$ year$^{-1}$ and $7.962 \times 10^5$ \$ year$^{-1}$, respectively. It is observed that many bypass efficiencies in the solution of $7.962 \times 10^5$ \$ year$^{-1}$ from Py-SLSQP are fractional, indicating premature termination of the optimisation. This is because distillation columns with fractional bypass efficiencies are thermodynamically inefficient (Dowling & Biegler, 2015). However, I-SLSQP needs longer computational time for more than 20 problem instances (i.e. half of the 40 instances), as shown in Fig. 4b. The difference between the areas on the two sides of the 0-level line in Fig. 4b is quite small, indicating that Py-SLSQP is slightly faster than I-SLSQP. The largest $|\gamma_t|$ appears in problem instance 1, which corresponds to solving DWCP from the first initial point. The total simulation time required from I-SLSQP and Py-SLSQP is

3890 s and 645 s, respectively.

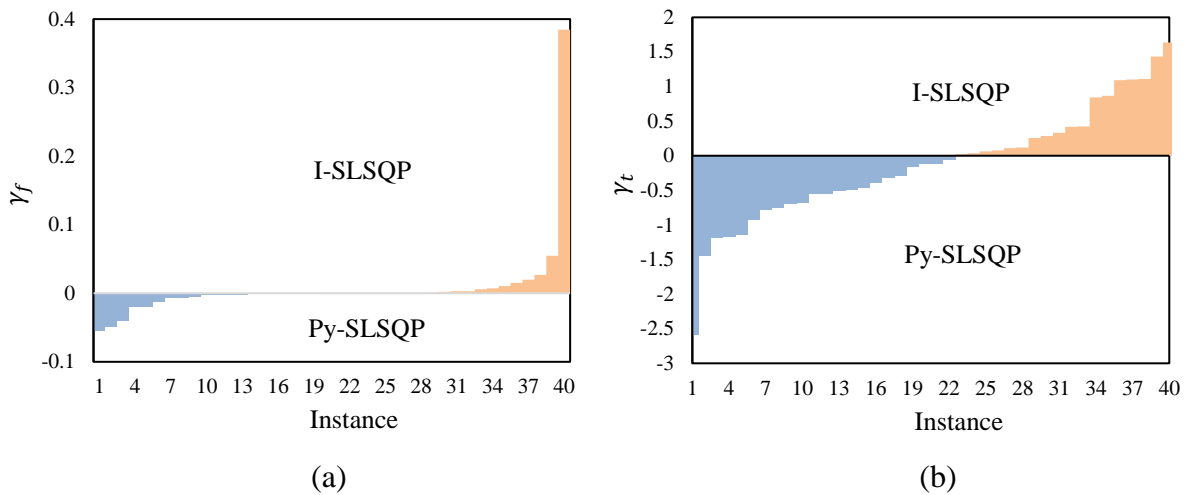

Figure 4 Morales profiles for 40 problem instances where (a) optima and (b) simulation time

In brief, I-SLSQP is more robust than Py-SLSQP as I-SLSQP is able to solve more problem instances. This is because I-SLSQP overcomes some limitations of the existing SLSQP algorithm in Py-SLSQP. In addition, I-SLSQP is less likely to terminate at a suboptimal point, although Py-SLSQP is faster.

## 5.3 I-SQP vs. I-SLSQP

As mentioned previously, SLSQP is expected to perform better than SQP when solving the ill-conditioned optimisation problems. To show this, we divide the optimisation problems into two sets in terms of the condition numbers of the reduced Hessian matrix $B_r$ at optima. The calculation of $B_r$ is shown in **Appendix B**. The problems PSD and DWC are included in set 1, whilst the other five problems are allocated to set 2 according to whether the average condition number $\bar{\kappa}$ of $Z^T B Z$ is smaller or bigger than $10^6$ as illustrated in Table 4. Here, the problems with $\bar{\kappa} \leq 10^6$ are considered to be well-conditioned; otherwise, the problems are ill-conditioned. The average number of function evaluations $\bar{n}_f$, and average simulation time $\bar{t}_{sim}$ for each problem are also provided in Table 4.

Table 4 Average condition numbers and computational performance for all seven problems

| Item | PSD | DWC | DWCP | DME | EDWC-AC | EDWC-EW | HPEDWC |
|---|---|---|---|---|---|---|---|
| $\bar{\kappa}$ | 3912 | $9.7 \times 10^4$ | $1.9 \times 10^7$ | $4.5 \times 10^8$ | $2.8 \times 10^{11}$ | $1.9 \times 10^6$ | $9.9 \times 10^{15}$ |
| $\bar{n}_f$ | 102 | 198 | 934 | 833 | 271 | 436 | 694 |
| $\bar{t}_{sim}$ | 120 | 201 | 1407 | 2367 | 1106 | 2058 | 2015 |

From Table 4, it can be observed that a much larger $\kappa$ ($> 10^6$) does exist for the last five problems. Therefore, it is much more difficult to solve these five optimisation problems, as indicated by both large values of $\bar{n}_f$ and $\bar{t}_{sim}$. For instance, the complexity of the optimisation model for the problem DWCP is very similar to that for DWC as shown in Table 1. However, $\bar{n}_f$ and $\bar{t}_{sim}$ for DWCP are around 5 and 7 times of those for DWC. The large difference in $\bar{n}_f$ and $\bar{t}_{sim}$ is mainly because the condition number increases by more than two orders of magnitude from DWC to DWCP.

### 5.3.1 I-SQP vs. I-SLSQP for well-conditioned problems

As seen from Fig. 5a, I-SLSQP and I-SQP generate very close solutions with a difference of less than 0.3% for the well-conditioned problems (i.e. the problems PSD, and DWC). The computational efficiency of I-SQP dominates that of I-SLSQP due to the much larger area in the half space below the $x$ axis in Fig. 5b. Specifically, I-SQP needs less simulation time for 8 out of 12 instances, especially for the first four instances in which the time savings are more than 50%.

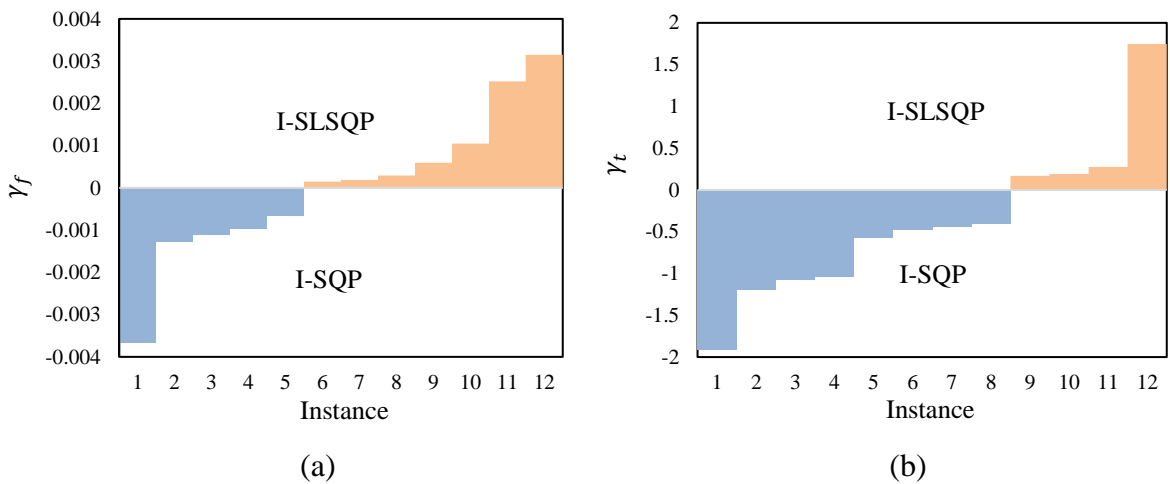

Figure 5 Morales profiles for 12 well-conditioned problem instances where (a) optima and (b) total simulation time.

### 5.3.2 I-SQP vs. I-SLSQP for ill-conditioned problems

As the algorithms perform evidently differently for the five ill-conditioned problems, we do not draw their Morales profiles together. Instead, we divide them into three groups: DWCP, DME and EDWC problems, and display the Morales profiles for each, such that the problem instances with similar profiles appear in the same figure. Fig. 6 illustrates the Morales profiles for DWCP problem instances.

As seen from Fig. 6a, I-SQP and I-SLSQP generate similar optimal solutions for most DWCP problem instances except for one, where I-SQP generates a bad local optimum with around 20% higher TAC. For computational efficiency, I-SQP is two to eight times faster than I-SLSQP for all six problem instances.

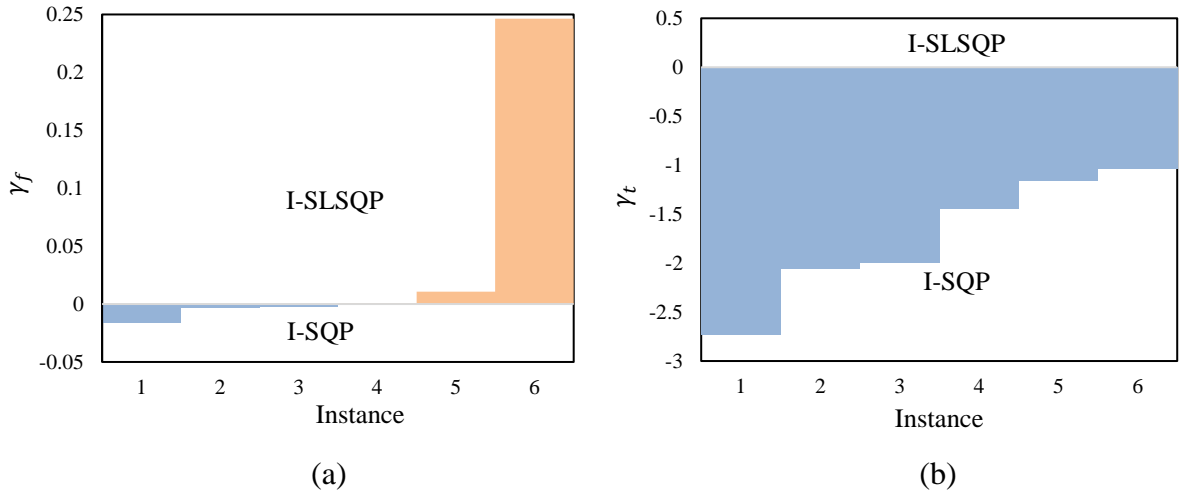

(a)            (b)

Figure 6 Morales profiles of I-SQP vs. I-SLSQP for six DWCP problem instances where (a) optima and (b) total simulation time

The Morales profiles for DME problem instances are illustrated in Fig. 7. From Fig. 7a, we can observe that I-SLSQP always generates better solutions, especially for the last two problem instances, where the optima from I-SQP are 40% and 100% higher than those from I-SLSQP. For computational efficiency, these two algorithms are basically competitive with each other according to Fig. 7b. In Fig. 8, we show the total simulation time for different problem instances which are in the same order as those in Fig. 7a. From Fig. 7a and Fig. 8, it is interestingly seen that I-SLSQP gets 50% lower TAC within around 70% less time for the 6$^{th}$ problem instance.

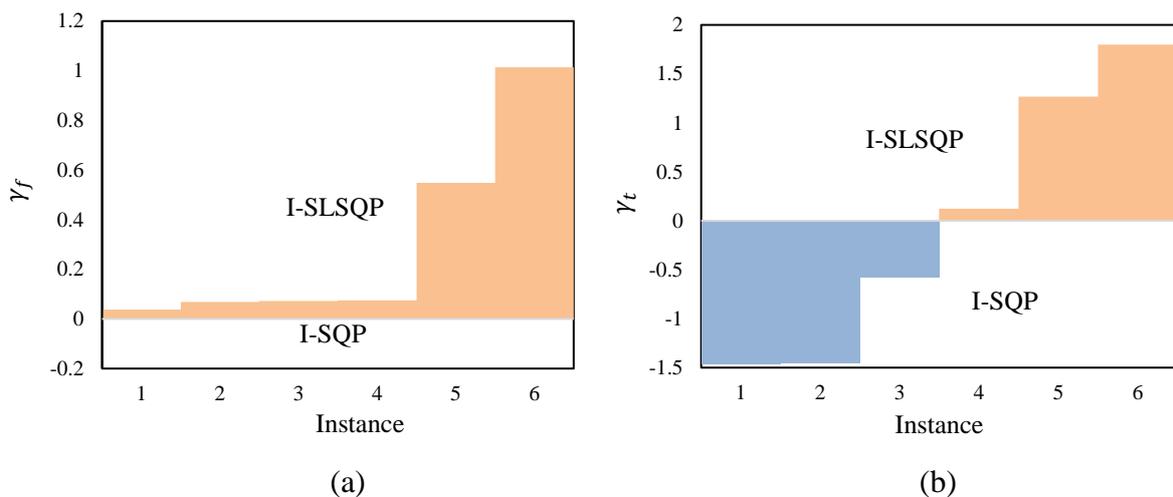

(a)            (b)

Figure 7. Morales profiles of I-SQP vs. I-SLSQP for six DME problem instances where (a) optima and (b) total simulation times

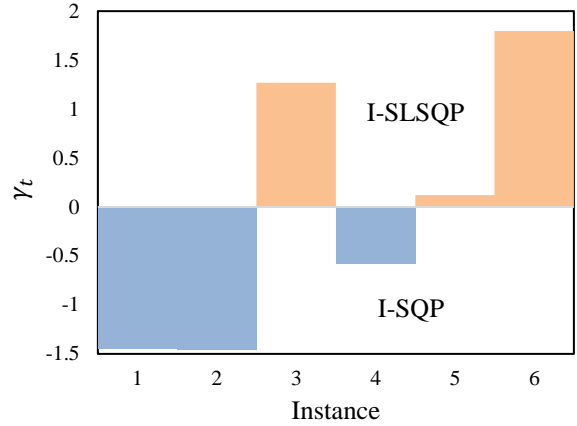

Figure 8 Logarithmic ratios of total simulation time between I-SQP and I-SLSQP in the order of ascending $\gamma_f$ for six DME problem instances

Finally, let's investigate the performance of these two algorithms when solving the three EDWC problems (i.e. EWDC-AC, EDWC-EW, and HPEDWC). As shown in Table 3, I-SQP fails to generate a solution for one EDWC problem instance due to the singular Jacobian matrix of the relaxed QP subproblem. As a result, the Morales profiles for only 17 EDWC problem instances are shown in Fig. 9. From Fig. 9a, it can be observed that I-SLSQP rarely generates worse solutions than I-SQP. From Fig. 9b, I-SLSQP needs less computational time in around two-thirds of the problem instances. The logarithmic ratios of total simulation time for the 17 problem instances in the order of ascending $\gamma_f$ are depicted in Fig. 10. From Fig. 10, I-SQP uses less computational time for instances 5, 7, and 9-12 where it terminates at a suboptimal point, while I-SLSQP achieves better solutions in less computational time for 7 problem instances (i.e. instances 6, 8, and 13-17).

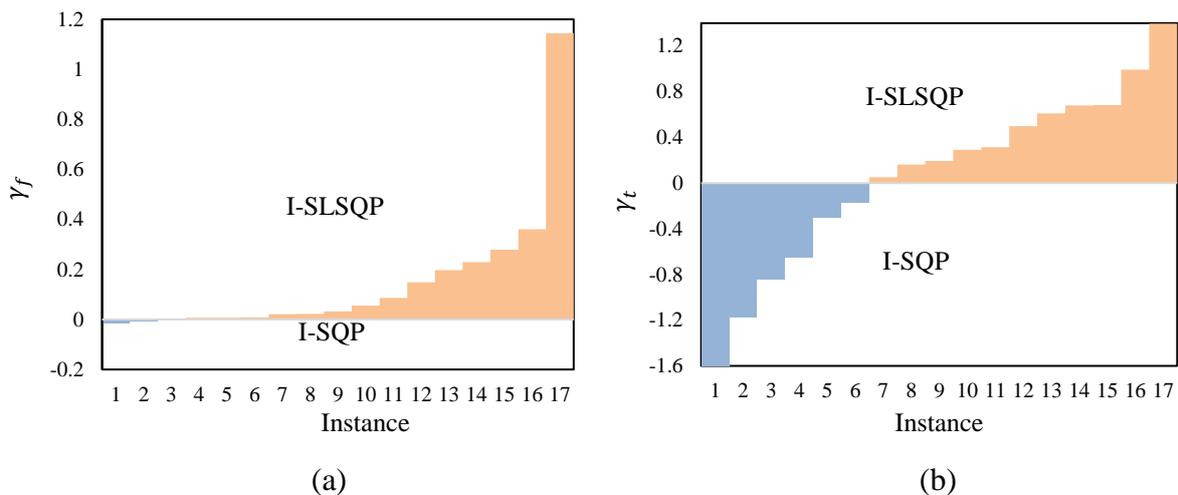

(a)          (b)

Figure 9 Morales profiles of I-SQP vs. I-SLSQP for 17 EDWC problem instances where (a) optima and (b) total simulation time

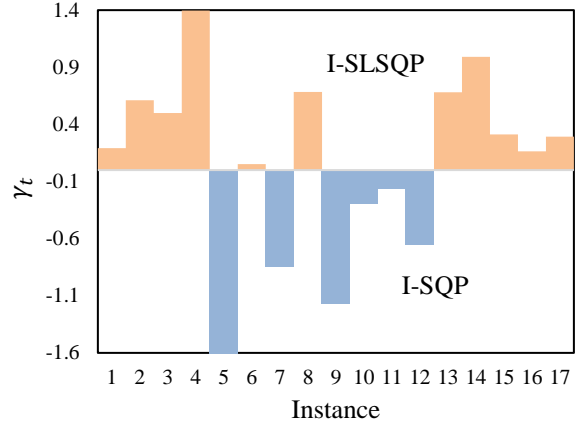

Figure 10 Logarithmic ratios of total simulation time between I-SQP and I-SLSQP for 17 EDWC problem instances in the order of ascending $\gamma_f$

## 5.4 Py-SLSQP vs. I-SQP

It is clearly seen that I-SQP often finds worse solutions for the last five ill-conditioned problems compared to I-SLSQP and Py-SLSQP, as shown in Tables 2 and 3. Therefore, there is no need to compare the performance of Py-SLSQP and I-SQP for those problems. We only make comparisons for the first two problems with 12 problem instances. The Morales profiles are shown in Fig. 11. As can be seen from Fig. 11, I-SQP is generally faster than Py-SLSQP, while the best solutions obtained from these two algorithms are very close.

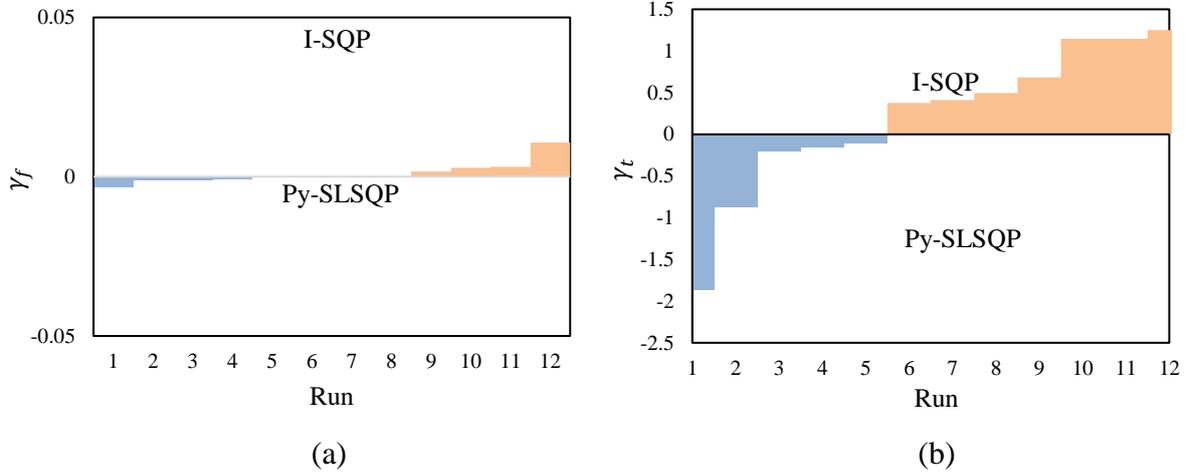

(a)  (b)

Figure 11 Morales profiles of Py-SLSQP vs. I-SQP for 12 problem instances of optimizing PSD, and DWC problems where (a) optima and (b) total simulation time.

## 5.5 Other discussions

We also observe the activation of the QP solver and the relaxations in I-SQP and I-SLSQP when solving these seven process optimisation problems from different initial points.

1. The activation of the QP solver in I-SLSQP. The QP solver is only activated when solving the EDWC-EW problem from the third and the fifth initial points. This is because the LSQ solver generates an ascent search direction even after the Hessian matrix is reset as identity. On the contrary, Py-SLSQP does not have such mechanism and hence terminates prematurely at an infeasible point when solving the EDWC-EW problem from the third initial point and when solving the HPEDWC example from the fourth initial point.

2. The usefulness of the relaxations of the QP/LSQ subproblems. In Table 5, we list the number of the original subproblems (QP/LSQ) and their relaxations solved when solving the seven large-scale process optimisation problems from different initial points. As seen from Table 5, the (RQP1) and (LSQ1) subproblems are required to generate optimal solutions for all problem instances, especially for the EDWC-EW and HPEDWC problems, where the relaxations (RQP1/LSQ1) are solved more than 10 times for most problem instances. We also observe that the relaxation LSQ2 is not activated at all when solving all the problems using I-SLSQP. However, the relaxation RQP2 is required when solving the EDWC-EW and HPEDWC problems using I-SQP. All these indicate that the use of the relaxations (RQP1/LSQ1) is able to address most of the inconsistent subproblems, whilst the relaxation RQP2/LSQ2 may be required in some instances.

3. I-SLSQP and Py-SLSQP have demonstrated better performance than I-SQP when they are used to solve ill-conditioned process optimisation problems. The possible reason is due to the use of the duality algorithm proposed by Lawson and Hanson (1995) for solving the LSQ subproblems in I-SLSQP and Py-SLSQP.

Table 5 Number of subproblems (QP/LSQ), (RQP1/LSQ1), and (RQP2/LSQ2) solved

| Problem | Algorithm | Initial points | | | | | |
|---|---|---|---|---|---|---|---|
| | | 1 | 2 | 3 | 4 | 5 | 6 |
| PSD | I-SQP | 44/4/0 | 122/1/0 | 56/1/0 | 61/2/0 | 71/2/0 | 65/2/0 |
| | I-SLSQP | 53/4/0 | 54/1/0 | 59/1/0 | 57/1/0 | 87/1/0 | 83/1/0 |
| DWC | I-SQP | 50/0/0 | 46/0/0 | 128/3/0 | 46/2/0 | 47/2/0 | 47/3/0 |
| | I-SLSQP | 61/0/0 | 77/0/0 | 92/1/0 | 51/2/0 | 62/2/0 | 65/3/0 |
| DWCP | I-SQP | 183/13/9 | 65/1/0 | 69/2/0 | 66/2/0 | 103/3/0 | 110/3/0 |
| | I-SLSQP | 349/0/0 | 74/0/0 | 122/1/0 | 202/1/0 | 115/1/0 | 88/1/0 |
| DME | I-SQP | 238/3/0 | 405/10/0 | 135/1/0 | 585/9/0 | 94/0/0 | 228/0/0 |
| | I-SLSQP | 188/0/0 | 309/0/0 | 319/1/0 | 421/0/0 | 168/0/0 | 303/0/0 |
| EDWC-AC | I-SQP | 208/6/0 | 95/2/0 | 123/2/0 | 184/4/0 | 130/9/0 | 222/15/0 |
| | I-SLSQP | 115/5/0 | 67/2/0 | 106/0/0 | 64/3/0 | 95/3/0 | 122/2/0 |
| EDWC-EW | I-SQP | 126/29/0 | 231/29/50 | 9/9/0 | 124/23/15 | 260/58/10 | 205/40/11 |
| | I-SLSQP | 114/7/0 | 189/19/0 | 82/41/0 | 80/25/0 | 136/73/0 | 193/31/0 |
| HPEDWC | I-SQP | 247/19/46 | 150/17/8 | 251/10/30 | 261/22/24 | 132/3/17 | 342/30/11 |
| | I-SLSQP | 256/14/0 | 173/24/0 | 331/23/0 | 311/32/0 | 152/10/0 | 187/41/0 |

Note: the three values divided by "/" in each cell correspond to the number of the subproblems (QP or LSQ) solved, the number of the relaxation (RQ1 or LSQ1) solved, and the number of the relaxation (RQ2 or LSQ2) solved, respectively.

## 6 Conclusion

In this work we developed an improved SQP algorithm, I-SQP and an improved SLSQP algorithm, I-SLSQP to drive the feasible path algorithm for process optimisation. To solve the inconsistent QP/LSQ subproblems reliably, a hybrid relaxation strategy was proposed, in which the modified Powell's relaxation formulation (RQP1)/(RLSQ1) with one relaxation variable was used in both algorithms at first, while the relaxation formulation with multiple relaxation variables (RQP2)/(RLSQ2) was activated if the former failed to generate a nonzero descent direction or the relaxation variable $\xi$ was close to 1 for more than certain continuous iterations. However, if the problem (RQP2)/(RLSQ2) led to an ill-conditioned Jacobian of the active constraints for more than certain continuous iterations, the formulation (RQP1)/(RLSQ1) would be switched back to in order to address the inconsistent problems. Our analysis of a part of the dual LSQ solution algorithm showed that it might suffer from serious cancellation errors, leading to the wrong solution or no solution when solving the LSQ subproblems in the SLSQP algorithm. In such cases, I-SLSQP would first reset $L^k$ and $D^k$ to the identity matrix and then resolve the LSQ subproblem. If the resetting still didn't work, the QP solver would be used to solve the LSQ subproblem.

Seven challenging process optimisation problems from the literature were solved from different initial points to illustrate the capabilities of I-SQP and I-SLSQP. The computational results showed that the existing Py-SLSQP algorithm and the proposed I-SQP and I-SLSQP algorithms had much better convergence performance than fmincon and IPOPT as the former three algorithms could solve most problems, while the latter two could only solve one or two out of the seven problems. It has also shown that I-SLSQP was more robust than I-SQP and Py-SLSQP as only I-SLSQP could solve all problem instances. I-SLSQP and the existing Py-SLSQP algorithm were competitive with each other in terms of solution quality and computational efficiency for problems that were solvable for both algorithms. I-SLSQP and Py-SLSQP usually generated much better solutions than I-SQP for ill-conditioned optimisation problems, especially for the EDWC problems, where I-SLSQP provided better solutions in shorter solution times for many instances. However, I-SQP was evidently faster than both SLSQP algorithms for well-conditioned problems.

**Acknowledgement**

The authors would like to appreciate the financial support from China Scholarship Council – The University of Manchester Joint Scholarship (201809120005), the President's Doctoral Scholar Award, The University of Manchester, and the PhD scholarship from Department of Chemical Engineering, The University of Manchester. Jie Li also thanks the financial support from UKRI Impact Acceleration Account (IAA 424).

# Appendix A Error analysis on the solution of LDP problem by dual algorithm

The solution of the problem (LDP) is crucial for solving problem (LSQ). Its numerical errors are analysed here.

The relative errors $\boldsymbol{\epsilon}^r$ of $\mathbf{r}$ got from Eq. (22) are bounded by Eq. A1 (Higham, 2002),

$$\boldsymbol{\epsilon}^r \leq diag(\max(\boldsymbol{\epsilon}^a, \boldsymbol{\epsilon}^b)) \frac{|\mathbf{a}^*| + |\mathbf{b}^*|}{|\mathbf{r}^*|}, \tag{A1}$$

where $\mathbf{a}^* = A^*\mathbf{u}^*$, and $\mathbf{r}^*$, $\mathbf{b}^*$, $A^*$ and $\mathbf{u}^*$ are the true values of $\mathbf{r}$, $\mathbf{b}$, $A$ and $\mathbf{u}$ when there is no numerical error. $diag(\cdot)$ is a function to construct a diagonal matrix from a given vector, while $\boldsymbol{\epsilon}^a \geq 0$ and $\boldsymbol{\epsilon}^b \geq 0$ are the relative errors of the actual $A\mathbf{u}$ and $\mathbf{b}$ respectively. In other words, $A\mathbf{u} = diag(A^*\mathbf{u}^*)(1 \pm \boldsymbol{\epsilon}^a)$ and $\mathbf{b} = \mathbf{b}^*(1 \pm \boldsymbol{\epsilon}^b)$. Note that $\mathbf{b} = [\overbrace{0, \ldots, 0}^{n'}, 1]^T$ is a constant vector, so $\boldsymbol{\epsilon}^b = 0$, $\mathbf{b}^* = \mathbf{b}$, and

$$r_i^* = a_i^*, i = 1, 2, \ldots, n', \tag{A2}$$

$$r_{n'+1}^* = a_{n'+1}^* - 1. \tag{A3}$$

From Eq. (A2), we have

$$\epsilon_i^r = \epsilon_i^a, i = 1, 2, \ldots, n', \tag{A4}$$

so the errors for the first $n'$ components of $\mathbf{r}$ are the same as the initial errors $\boldsymbol{\epsilon}^a$.

Substitute Eq. (A3) to Eq. (A1), we get

$$\epsilon_{n'+1}^r \leq \epsilon_{n'+1}^a \frac{|a_{n'+1}^*| + 1}{r_{n'+1}^*} = \epsilon_{n'+1}^a \frac{|1 - r_{n'+1}^*| + 1}{r_{n'+1}^*}. \tag{A5}$$

Eq. (A5) can be simplified as

$$\epsilon_{n'+1}^r \leq \epsilon_{n'+1}^a \left(\frac{2}{r_{n'+1}^*} - 1\right) \text{ when } r_{n'+1}^* < 1. \tag{A6}$$

Eq. (A6) can be further simplified as

$$\epsilon_{n'+1}^r \leq 2 \frac{\epsilon_{n'+1}^a}{r_{n'+1}^*} \text{ when } r_{n'+1}^* \ll 1. \tag{A7}$$

Eq. (A7) indicates that $\epsilon_{n'+1}^r$ may be significantly larger than $\epsilon_{n'+1}^a$ when $r_{n'+1}^* \ll 1$. Note that when $r_{n'+1} \leq 0$, the problem (LDP) is declared as infeasible (Lawson and Hanson, 1995), and neither do the problem (LSQ). Since $r_{n'+1} = r_{n'+1}^*(1 \pm \epsilon_{n'+1}^r)$, it may be nonpositive when $r_{n'+1}^* \leq 2\epsilon_{n'+1}^1$ and $\epsilon_{n'+1}^r \geq 1$ even if $r_{n'+1}^* > 0$.

From Eq. (23), the errors in $\mathbf{r}$ will propagate into $\mathbf{z}$ as shown in the following Eq. (A8).

$$z_i = -\frac{r_i^*(1 \pm \epsilon_i^a)}{r_{n'+1}^*\left(1 \pm 2\frac{\epsilon_{n'+1}^a}{r_{n'+1}^*}\right)} = z_i^* \frac{(1 \pm \epsilon_i^a)}{\left(1 \pm 2\frac{\epsilon_{n'+1}^a}{r_{n'+1}^*}\right)}, i = 1, 2, \ldots, n' \text{ when } r_{n'+1}^* \ll 1, \tag{A8}$$

where $z_i^* = -\dfrac{r_i^*}{r_{n'+1}^*}, i = 1, 2, \ldots, n'$.

## Appendix B Calculation of the reduced Hessian matrix

The reduced Hessian matrix $B_r$ is got from

$$B_r = Z^T B Z, \tag{B1}$$

where $Z \in \mathbb{R}^{n \times (n-m)}$ is the null space matrix of the Jacobian matrix for active constraint $j \in \mathcal{A}$, i.e.,

$$\bar{A} = \left[\nabla c_j\right]_{j \in \mathcal{A}}^T \tag{B2}$$

$$\bar{A} Z = 0. \tag{B3}$$

Here, $\mathbf{c} = [\mathbf{g}^T, \mathbf{h}^T]^T$ and $\bar{A}$ is the Jacobian matrix for active constraints.

## Appendix C Initial values of decision variables for all seven process optimisation problems

Table C1 Initial values of decision variables for the problem PSD

| Variable | Value | Variable | Value |
| --- | --- | --- | --- |
| $A_1$ (m²) | 100.0 | $L_1$ (kmol h⁻¹) | 1880.0 |
| $A_2$ (m²) | 100.0 | $VF$ (kmol kmol⁻¹) | 0.7 |
| $P_{LP}$ (bar) | 1.1 | $RR_1$ (kmol kmol⁻¹) | 1.0 |
| $P_{HP}$ (bar) | 10.0 | $RR_2$ (kmol kmol⁻¹) | 1.5 |
| $\epsilon_j, j = 1,2, \cdots 54$ | 0.1, 0.3, 0.5, 0.7, 0.9, 1.0 | | |

$A_1$: area of heat exchanger 1; $A_2$: area of heat exchanger 2; $P_{LP}$: pressure of the low pressure (LP) column; $P_{HP}$: pressure of the high pressure (HP) column; $L_1$: bottom flow rate of the LP column; VF: reboiler vaporisation fraction (VF); $RR_1$: reflux ratio of the LP column; $RR_2$: reflux ratio of the HP column; $\epsilon_j$: bypass efficiency at stage $j$.

Table C2 Initial values of decision variables for the problem DWC

| Variable | Value | Variable | Value |
| --- | --- | --- | --- |
| $D$ (kmol h⁻¹) | 150 | $SF_1$ (kmol kmol⁻¹) | 0.5 |
| $VF$ (kmol kmol⁻¹) | 0.6 | $SF_2$ (kmol kmol⁻¹) | 0.5 |
| $SD$ (kmol kmol⁻¹) | 0.5 | | |
| $\epsilon_j, j = 1,2, \cdots 180$ | 0.1, 0.3, 0.5, 0.7, 0.9, 1.0 | | |

$D$: distillate flow rate; $VF$: reboiler vaporisation fraction; $SD$: side draw fraction; $SF_1$: liquid split fraction from the column top to the left-hand side of the dividing wall; $SF_2$: vapour split fraction from the column bottom to the right-hand side of the dividing wall; $\epsilon_j$: bypass efficiency at stage $j$.

Table C3 Initial values of decision variables for the problem DWCP

| Variable | Value | Variable | Value |
|---|---|---|---|
| $D$ (kmol h$^{-1}$) | 150 | $SF_1$ (kmol kmol$^{-1}$) | 0.5 |
| $VF$ (kmol kmol$^{-1}$) | 0.6 | $SF_2$ (kmol kmol$^{-1}$) | 0.5 |
| $SD$ (kmol kmol$^{-1}$) | 0.5 | $P$ (bar) | 1.2 |
| $\epsilon_j, j = 1,2,\cdots 180$ | 0.1, 0.3, 0.5, 0.7, 0.9, 1.0 | | |

$D$: flow rate; $VF$: reboiler vaporisation fraction; $SD$: side draw fraction; $SF_1$: liquid split fraction from the column top to the left-hand side of the dividing wall; $SF_2$: vapour split fraction from the column bottom to the right-hand side of the dividing wall; $P$: column pressure; $\epsilon_j$: bypass efficiency at stage $j$.

Table C4 Initial values of decision variables for the problem DME

| Variable | Value | Variable | Value |
|---|---|---|---|
| $F$ (kmol h$^{-1}$) | 920 | $P_{COL}$ (bar) | |
| $VFV$ (kmol kmol$^{-1}$) | 1.1 | $RR$ (kmol kmol$^{-1}$) | |
| $T_{HX}$ (°C) | 300 | $VF$ (kmol kmol$^{-1}$) | 0.6 |
| $V_R$ (m$^3$) | 100 | $SD$ (kmol kmol-1) | 0.5 |
| $P_R$ (bar) | 10 | $SF_1$ (kmol kmol$^{-1}$) | 0.5 |
| $T_{CO}$ (°C) | 150 | $SF_2$ (kmol kmol$^{-1}$) | 0.5 |
| $\epsilon_j, j = 1,2,\cdots 110$ | 0.1, 0.3, 0.5, 0.7, 0.9, 1.0 | | |

$F$: flow rate of fresh methanol; $VFV$: vaporisation fraction in the vaporiser; $T_{HX}$: outlet temperature of the cold stream in the heat exchanger; $V_R$: volume of the reactor; $P_R$: operating pressure of the distillation column; $T_{CO}$: temperature of the cooler; $P_{COL}$: operating pressure of the distillation column; $RR$: reflux ratio; $VF$: vaporisation fraction in the reboiler; $SD$: side draw fraction; $SF_1$: liquid split fraction from the column top section to the left-hand side of the dividing wall; $SF_2$: vapour split fraction from the column bottom section to the right-hand side of the dividing wall; $\epsilon_j$: bypass efficiency at stage $j$.

Table C5 Initial values of decision variables for the problem EDWC-AW

| Variable | Value | Variable | Value |
|---|---|---|---|
| $F_E$ (kmol h$^{-1}$) | 0.01 | $SF$ (kmol kmol-1) | 0.5 |
| $RR_M$ (kmol kmol$^{-1}$) | 1.0 | $F_B$ (kmol h$^{-1}$) | 100 |
| $RR_S$ (kmol kmol$^{-1}$) | 1.0 | | |
| $\epsilon_j, j = 1,2,\cdots 85$ | 0.1, 0.3, 0.5, 0.7, 0.9, 1.0 | | |

$F_E$: the entrainer make-up flow rate; $RR_M$: reflux ratio of the main column; $RR_S$: reflux ratio of the side column; $SF$: split fraction of the vapour stream to side column; $F_B$: column bottom flow rate; $\epsilon_j$: bypass efficiency at stage $j$.

Table C6 Initial values of decision variables for the problem EDWC-EW

| Variable | Value | Variable | Value |
|---|---|---|---|
| $F_E$ (kmol h$^{-1}$) | 0.01 | $SF$ (kmol kmol-1) | 0.5 |
| $RR_M$ (kmol kmol$^{-1}$) | 1.0 | $F_B$ (kmol h$^{-1}$) | 300 |
| $VF$ (kmol kmol$^{-1}$) | 0.5 | | |
| $\epsilon_j, j = 1,2,\cdots 130$ | 0.1, 0.3, 0.5, 0.7, 0.9, 1.0 | | |

$F_E$: the entrainer make-up flow rate; $RR_M$: reflux ratio of the main column; $VF$: vaporisation fraction in the reboiler; $SF$: split fraction of the vapour stream to side column; $F_B$: column bottom flow rate; $\epsilon_j$: bypass efficiency at stage $j$.

Table C7 Initial values of decision variables for the problem HPEDWC

| Variable | Value | Variable | Value |
|---|---|---|---|
| $F_E$ (kmol h$^{-1}$) | 0.01 | $F_B$ (kmol h$^{-1}$) | 300 |
| $RR_M$ (kmol kmol$^{-1}$) | 1.0 | $P$ (atm) | 3 |
| $VF$ (kmol kmol$^{-1}$) | 0.5 | $A_1$ | 100 |
| $SF$ (kmol kmol-1) | 0.5 | $A_2$ | 10 |
| $\epsilon_j, j = 1,2,\cdots 130$ | 0.1, 0.3, 0.5, 0.7, 0.9, 1.0 | | |

$F_E$: the entrainer make-up flow rate; $RR_M$: reflux ratio of the main column; $VF$: vaporisation fraction in the reboiler; $SF$: split fraction of the vapour stream to side column; $F_B$: column bottom flow rate; $P$: outlet pressure of the compressor; $A_1$: area of the heat exchanger 1; $A_2$: area of the heat exchanger 2; $\epsilon_j$: bypass efficiency at stage $j$.

**References**

Amestoy, P.R., Duff, I.S., l'Excellent, J.-Y., 2000. Multifrontal parallel distributed symmetric and unsymmetric solvers. Computer Methods in Applied Mechanics and Engineering. 184 (2-4), 501-520.

Aspen Technology Inc., 2015. Aspen Custom Modeler User's Guide, in: Technology, A. (Ed.), http://www.aspentech.com.

Biegler, L.T., 1993. From nonlinear programming theory to practical optimization algorithms: A process engineering viewpoint. Comput. Chem. Eng. 17, S63-S80.

Biegler, L.T., 2010. Nonlinear programming: concepts, algorithms, and applications to chemical processes. Society for Industrial and Applied Mathematics, Philadelphia, Pennsylvania.

Biegler, L.T., Cuthrell, J.E., 1985. Improved infeasible path optimization for sequential modular simulators—II: the optimization algorithm. Comput. Chem. Eng. 9, 257-267.

Biegler, L.T., Hughes, R.R., 1981. Approximation programming of chemical processes with Q/LAP. Chemical Engineering Progress 77.

Biegler, L.T., Hughes, R.R., 1982. Infeasible path optimization with sequential modular simulators. AIChE J. 28, 994-1002.

Biegler, L.T., Hughes, R.R., 1985. Feasible path optimization with sequential modular simulators. Comput. Chem. Eng. 9, 379-394.

Boggs, P.T., Tolle, J.W., 1995. Sequential quadratic programming. Acta Numerica 4, 1-51.

Byrd, R.H., Hribar, M.E., Nocedal, J., 1999. An Interior Point Algorithm for Large-Scale Nonlinear Programming. SIAM Journal on Optimization 9, 877-900.

Caballero, J.A., Grossmann, I.E., 2008. An algorithm for the use of surrogate models in modular flowsheet optimization. AIChE Journal 54, 2633-2650.

Dai, Y.-H., Schittkowski, K., 2008. A sequential quadratic programming algorithm with non-monotone line search. Pacific Journal of Optimization 4, 335-351.

Dowling A.W., Biegler, L.T., 2015. Rigorous Optimization-based Synthesis of Distillation Cascades without Integer Variables. In: Jiˇ ríJaromír Klemeš, P.S.V., Peng Yen, L. (Eds.), Computer Aided Chemical Engineering. Elsevier, pp. 55–60 .

Drud, A.S., 1994. CONOPT—A Large-Scale GRG Code. ORSA Journal on Computing 6, 207-216.

Fletcher, R., Powell, M.J.D., 1974. On the modification of $LDL^T$ factorizations. Mathematics of Computation 28, 1067-1087.

Gill, P.E., Hammarling, S.J., Murray, W., Saunders, M.A., Wright, M.H., 1986. User's Guide

for LSSOL (Version 1. 0): Fortran package for constrained linear least-squares and convex quadratic programming. Systems Optimization Lab., Stanford Univ., CA (USA).

Gill, P.E., Murray, W., Saunders, M.A., 2002. SNOPT: An SQP Algorithm for Large-Scale Constrained Optimization. SIAM Journal on Optimization 12, 979-1006.

Gill, P.E., Murray, W., Wright, M.H., 2019. Practical optimization. Society for Industrial and Applied Mathematics SIAM, 3600 Market Street, Philadelphia, PA 19104, Philadelphia, Pennsylvania.

Golub, G.H., Van Loan, C.F., 2013. Matrix computations, Fourth edition. ed. Johns Hopkins University Press, Baltimore, Maryland.

Gurobi Optimization, LLC, 2022. Gurobi optimizer reference manual, https://www.gurobi.com.

Han, S.-P., 1976. Superlinearly convergent variable metric algorithms for general nonlinear programming problems. Mathematical Programming 11, 263-282.

Han, S.P., 1977. A globally convergent method for nonlinear programming. Journal of Optimization Theory and Applications 22, 297-309.

Haskell, K.H., Hanson, R.J., 1981. An algorithm for linear least squares problems with equality and nonnegativity constraints. Mathematical Programming 21, 98-118.

Higham, N.J., 2002. Accuracy and stability of numerical algorithms, 2nd ed. ed. Society for Industrial and Applied Mathematics SIAM, 3600 Market Street, Floor 6, Philadelphia, PA 19104, Philadelphia, Pa.

Kossack, S., Kraemer, K., Marquardt, W., 2006. Efficient Optimization-Based Design of Distillation Columns for Homogenous Azeotropic Mixtures. Industrial & Engineering Chemistry Research 45, 8492-8502.

Kraft, D., 1988. A software package for sequential quadratic programming. Report DFVLR-FR 88–28 (Deutsche Forschungs- und Versuchsanstalt für Luftund Raumfahrt).

Lawson, C.L., Hanson, R.J., 1995. Solving least squares problems. SIAM, Philadelphia.

Ledezma-Martínez, M., Jobson, M., Smith, R., 2018. Simulation–Optimization-Based Design of Crude Oil Distillation Systems with Preflash Units. Industrial & Engineering Chemistry Research 57, 9821-9830.

Ma, Y., Luo, Y., Yuan, X., 2019. Towards the really optimal design of distillation systems: Simultaneous pressures optimization of distillation systems based on rigorous models. Comput. Chem. Eng. 126, 54-67.

Ma, Y., McLaughlan, M., Zhang, N., Li, J., 2020a. Novel feasible path optimisation algorithms using steady-state and/or pseudo-transient simulations. Comput. Chem. Eng. 143, 107058.

Ma, Y., Yang, Z., El-Khoruy, A., Zhang, N., Li, J., Zhang, B., Sun, L., 2021. Simultaneous

Synthesis and Design of Reaction–Separation–Recycle Processes Using Rigorous Models. Industrial & Engineering Chemistry Research 60, 7275-7290.

Ma, Y., Zhang, N., Li, J., Cao, C., 2020b. Optimal design of extractive dividing-wall column using an efficient equation-oriented approach. Front. Chem. Sci. Eng. 15, 72-89.

Morales, J.L., Nocedal, J., 2011. Remark on "algorithm 778: L-BFGS-B: Fortran subroutines for large-scale bound constrained optimization". ACM Transactions on Mathematical Software 38 (1), 1-4.

Moré, J.J., Wild, S.M., 2011. Estimating computational noise. SIAM Journal on Scientific Computing 33, 1292-1314.

Nocedal, J., Wright, S.J., 2006. Numerical optimization, Second Edition. ed. Springer New York, New York, NY.

Nowak, I., 1988. Ein quadratisches Optimierungsproblem mit schlupfvariablen fiir die SQPMethode zur l6sung des allgemeinen nichtlinearen optimierungsproblems. TH Darmstadt.

Oztoprak, F., Byrd, R., Nocedal, J.J.a.e.-p., 2021. Constrained Optimization in the Presence of Noise, p. arXiv:2110.04355.

Pattison, R.C., Baldea, M., 2014. Equation-oriented flowsheet simulation and optimization using pseudo-transient models. AIChE Journal 60, 4104-4123.

Gill, P., Murray, W., Saunders, M.., Wright, M.H., 1984. User's guide for OPSOL (Version 3.2): A Fortran package for quadratic programming.

Powell, M.J.D., 1978a. The convergence of variable metric methods for nonlinearly constrained optimization calculations, in: Mangasarian, O.L., Meyer, R.R., Robinson, S.M. (Eds.), Nonlinear Programming 3. Academic Press, pp. 27-63.

Powell, M.J.D., 1978b. A fast algorithm for nonlinearly constrained optimization calculations. Springer Berlin Heidelberg, Berlin, Heidelberg, pp. 144-157.

Powell, M.J.D.J.M.P., 1978c. Algorithms for nonlinear constraints that use lagrangian functions. Mathematical Programming 14, 224-248.

Python Software Foundation, 2016. Python language reference Version 3.6.

Schittkowski, K., 1980. Nonlinear programming codes. Springer-Verlag, Heidelberg.

Schittkowski, K., 1982. The nonlinear programming method of Wilson, Han, and Powell with an augmented Lagrangian type line search function. Numerische Mathematik 38, 115-127.

Schittkowski, K., 2008. An updated set of 306 test problems for nonlinear programming with validated optimal solutions—user's guide. Research report, Department of Computer Science, University of Bayreuth.

Schittkowski, K., 2011. A robust implementation of a sequential quadratic programming algorithm with successive error restoration. Optimization Letters 5, 283-296.

Stoer, J., 1971. On the numerical solution of constrained least-squares problems. SIAM Journal on Numerical Analysis 8, 382-411.

The Mathworks, I., 2023. Matlab R2023b. The Mathworks, Inc., Natick, Massachusetts.

Tone, K., 1983. Revisions of constraint approximations in the successive QP method for nonlinear programming problems. Mathematical Programming 26, 144-152.

Virtanen, P., Gommers, R., Oliphant, T.E., Haberland, M., Reddy, T., Cournapeau, D., Burovski, E., Peterson, P., Weckesser, W., Bright, J., van der Walt, S.J., Brett, M., Wilson, J., Millman, K.J., Mayorov, N., Nelson, A.R.J., Jones, E., Kern, R., Larson, E., Carey, C.J., Polat, İ., Feng, Y., Moore, E.W., VanderPlas, J., Laxalde, D., Perktold, J., Cimrman, R., Henriksen, I., Quintero, E.A., Harris, C.R., Archibald, A.M., Ribeiro, A.H., Pedregosa, F., van Mulbregt, P., Vijaykumar, A., Bardelli, A.P., Rothberg, A., Hilboll, A., Kloeckner, A., Scopatz, A., Lee, A., Rokem, A., Woods, C.N., Fulton, C., Masson, C., Häggström, C., Fitzgerald, C., Nicholson, D.A., Hagen, D.R., Pasechnik, D.V., Olivetti, E., Martin, E., Wieser, E., Silva, F., Lenders, F., Wilhelm, F., Young, G., Price, G.A., Ingold, G.-L., Allen, G.E., Lee, G.R., Audren, H., Probst, I., Dietrich, J.P., Silterra, J., Webber, J.T., Slavič, J., Nothman, J., Buchner, J., Kulick, J., Schönberger, J.L., de Miranda Cardoso, J.V., Reimer, J., Harrington, J., Rodríguez, J.L.C., Nunez-Iglesias, J., Kuczynski, J., Tritz, K., Thoma, M., Newville, M., Kümmerer, M., Bolingbroke, M., Tartre, M., Pak, M., Smith, N.J., Nowaczyk, N., Shebanov, N., Pavlyk, O., Brodtkorb, P.A., Lee, P., McGibbon, R.T., Feldbauer, R., Lewis, S., Tygier, S., Sievert, S., Vigna, S., Peterson, S., More, S., Pudlik, T., Oshima, T., Pingel, T.J., Robitaille, T.P., Spura, T., Jones, T.R., Cera, T., Leslie, T., Zito, T., Krauss, T., Upadhyay, U., Halchenko, Y.O., Vázquez-Baeza, Y., SciPy, C., 2020. SciPy 1.0: fundamental algorithms for scientific computing in Python. Nature Methods 17, 261-272.

A. Wächter and L. T. Biegler, 2006. On the Implementation of a Primal-Dual Interior Point Filter Line Search Algorithm for Large-Scale Nonlinear Programming, Mathematical Programming 106(1), 25-57.

Wilson, R.B., 1963. A simplicial algorithm for concave programming. Graduate School of Business Administration, George F. Baker Foundation, Harvard University, Boston, Mass.